\documentclass{amsart}

\usepackage{amsmath,amssymb,url,color,mathtools}
\usepackage[margin=4cm]{geometry}
\usepackage{tikz}

\newtheorem{lemma}{Lemma}[section]
\newtheorem{prop}[lemma]{Proposition}
\newtheorem{cor}[lemma]{Corollary}
\newtheorem{thm}[lemma]{Theorem}
\newtheorem{example}[lemma]{Example}
\newtheorem{thm?}[lemma]{Theorem?}

\newtheorem{remark}[lemma]{Remark}

\begin{document}
\title[The Set of Functonal Degrees]{Functional Degrees and Arithmetic Applications, I: The Set of Functonal Degrees}

\author{Pete L. Clark}
\author{Uwe Schauz}

\theoremstyle{definition}
\newtheorem{dfn}[lemma]{Definition}
\renewcommand\atop[2]{\genfrac{}{}{0pt}{}{#1}{#2}}
\newcommand{\Mod}[1]{\ (\mathrm{mod}\ #1)}
\newcommand{\etalchar}[1]{$^{#1}$}
\newcommand{\F}{\mathbb{F}}
\newcommand{\et}{\textrm{\'et}}
\newcommand{\ra}{\ensuremath{\rightarrow}}
\newcommand{\lra}{\ensuremath{\longrightarrow}}
\newcommand{\FF}{\F}
\newcommand{\ff}{\mathfrak{f}}
\newcommand{\Z}{\mathbb{Z}}
\newcommand{\N}{\mathbb{N}}
\newcommand{\NN}{\widetilde{\N}}
\newcommand{\mm}{\underline{m}}
\newcommand{\nn}{\underline{n}}
\newcommand{\ch}{}
\newcommand{\R}{\mathbb{R}}
\renewcommand{\P}{\mathbb{P}}
\newcommand{\PP}{\mathbb{P}}
\newcommand{\pp}{\mathfrak{p}}
\newcommand{\C}{\mathbb{C}}
\newcommand{\Q}{\mathbb{Q}}
\newcommand{\ab}{\operatorname{ab}}
\newcommand{\Aut}{\operatorname{Aut}}
\newcommand{\gk}{\mathfrak{g}_K}
\newcommand{\gq}{\mathfrak{g}_{\Q}}
\newcommand{\OQ}{\overline{\Q}}
\newcommand{\Out}{\operatorname{Out}}
\newcommand{\End}{\operatorname{End}}
\newcommand{\Gal}{\operatorname{Gal}}
\newcommand{\CT}{(\mathcal{C},\mathcal{T})}
\newcommand{\lcm}{\operatorname{lcm}}
\newcommand{\Div}{\operatorname{Div}}
\newcommand{\OO}{\mathcal{O}}
\newcommand{\rank}{\operatorname{rank}}
\newcommand{\tors}{\operatorname{tors}}
\newcommand{\IM}{\operatorname{IM}}
\newcommand{\CM}{\mathbf{CM}}
\newcommand{\HS}{\mathbf{HS}}
\newcommand{\Frac}{\operatorname{Frac}}
\newcommand{\Pic}{\operatorname{Pic}}
\newcommand{\coker}{\operatorname{coker}}
\newcommand{\Cl}{\operatorname{Cl}}
\newcommand{\loc}{\operatorname{loc}}
\newcommand{\GL}{\operatorname{GL}}
\newcommand{\PGL}{\operatorname{PGL}}
\newcommand{\PSL}{\operatorname{PSL}}
\newcommand{\Frob}{\operatorname{Frob}}
\newcommand{\Hom}{\operatorname{Hom}}
\newcommand{\Coker}{\operatorname{\coker}}
\newcommand{\Ker}{\ker}
\newcommand{\g}{\mathfrak{g}}
\newcommand{\sep}{\operatorname{sep}}
\newcommand{\new}{\operatorname{new}}
\newcommand{\Ok}{\mathcal{O}_K}
\newcommand{\ord}{\operatorname{ord}}
\newcommand{\Ohell}{\OO_{\ell^{\infty}}}
\newcommand{\cc}{\mathfrak{c}}
\newcommand{\ann}{\operatorname{ann}}
\renewcommand{\tt}{\mathfrak{t}}
\renewcommand{\cc}{\mathfrak{a}}
\renewcommand{\aa}{\mathfrak{a}}
\newcommand\leg{\genfrac(){.4pt}{}}
\renewcommand{\gg}{\mathfrak{g}}
\renewcommand{\O}{\mathcal{O}}
\newcommand{\Spec}{\operatorname{Spec}}
\newcommand{\rr}{\mathfrak{r}}
\newcommand{\rad}{\operatorname{rad}}
\newcommand{\SL}{\operatorname{SL}}
\newcommand{\fdeg}{\operatorname{fdeg}}
\renewcommand{\rank}{\operatorname{rank}}
\newcommand{\Int}{\operatorname{Int}}
\newcommand{\zz}{\mathbf{z}}

\begin{abstract}
We give a further development of the Aichinger-Moosbauer calculus of functional
degrees of maps between commutative groups.
%which is based on the simple idea that the functional degree should decrease
%if a discrete derivative is taken.
For any fixed given commutative groups $A$ and $B$, we compute the largest
possible finite functional degree that a map $f: A \ra B$ can have. We also
determine the set of all possible degrees of such maps. This also yields a
solution to Aichinger and Moosbauer's problem of finding the nilpotency index of
the augmentation ideal of group rings of the form
$Z_{p^\beta}[Z_{p^{\alpha_1}}\!\times Z_{p^{\alpha_2}} \!\times\dotsm\times
Z_{p^{\alpha_n}}]$ with $p,\beta,n,\alpha_1,\dotsc,\alpha_n\in\Z^+$\!, $p$ prime.
%Next, for any commutative group $B$,
%we give a canonical series representation for any function $f: \Z^N \ra B$ of finite
%functional degree.  Finally, we combine the aforementioned results with work of Wilson
%\cite{Wilson06} to give a group-theoretic result that recovers the theorems on higher
%$p$-adic congruences for the size of the solution set of a polynomial system over
%$\F_q$ due to Ax-Katz and Moreno-Moreno.  Our work extends these theorems to any
%finite commutative ring of prime characteristic.
\end{abstract}

\maketitle

\tableofcontents

\section{Introduction}
\subsection{Polynomial Maps} This is the first of two papers in which we attempt a synthesis and
development of two prior works: a 2006 paper of R.\ Wilson \cite{Wilson06} and a
2021 paper of E.\ Aichinger and J.\ Moosbauer \cite{Aichinger-Moosbauer21}.
These works revisit and generalize number-theoretic results in connection with
polynomials, using methods that seem to belong to pure algebra.
\\ \\
The major algebraic theme is the notion of a ``polynomial map'' $f: A \ra B$ between
commutative groups $A$ and $B$.  It is an important concept that arises in many
contexts, for instance in the notion of a quadratic map between commutative groups: see
e.g.\ \cite{Zarkhin74} and \cite[p.\,34-35]{Serre}.  Over commutative groups $A$ and
$B$, however, there is not a set standard definition for ``polynomial maps and their
degree,'' but various different definitions. Indeed, many concurrent definitions are
compared in a work of Laczkovich \cite{Laczkovich04}.
%This approach was taken by Wilson in \cite{Wilson06} and by Aichinger and Moosbauer in
%\cite{Aichinger-Moosbauer21}.  It is the one of interest to us here.
\\ \\
One frequently used, simple and successful way of viewing polynomial maps comes
from the calculus of finite differences.  In 1909 Fr\'echet showed \cite{Frechet09} that
a continuous function $f: \R \ra \R$ is a polynomial if and only if for some $d \in \N$
we have $\Delta^d f = 0$, where $\Delta$ is the forward \textbf{difference operator}
with the defining equation
\[ (\Delta f)(x) = f(x+1)-f(x)\]
and $\Delta^d$ is its $d$-fold iterate.  Moreover, for a nonzero polynomial $f$ over $\R$,
the least $d$ such that $\Delta^d f = 0$ is equal to $\deg(f)+1$.  It is this property
that can be used to generalize the degree of polynomials to more general functions
between arbitrary commutative groups.  One simply has to extend the calculus of
finite differences to such functions. Various authors have taken on aspects of this
task.  A particularly systematic and penetrating take was recently given in work
 of Aichinger and Moosbauer \cite{Aichinger-Moosbauer21}.  Using finite differences
over groups, the authors associate, to any map $f: A \ra B$, a \textbf{functional
degree} $\fdeg(f)$ that is either a non-negative integer or $\infty$. They simply set
$$\fdeg(f):=\min\{d\in\Z^+\mid \Delta^d f = 0\}-1\,,$$ a definition that we follow,
except that, for the zero function, we set
$$\fdeg(0):=-\infty\,.$$
%\footnote{We also set $\fdeg(0):=-\infty$, in contrast to Achinger and Moosbauer's
%convention of having degree $0$ for all constant functions.}.
%For a field $F$ and a polynomial map $E(P): F^n \ra F$ -- i.e., a map obtained by evaluating a polynomial
%$P \in F[t_1,\ldots,t_n]$ with coefficients in $F$ -- they compare the functional degree of $E(P)$ to the degree of $P$ in the usual algebraic sense, showing in particular that equality always holds if and only if $F$ has characteristic $0$.

\subsection{Main Results}
Basic properties of the functional degree, such as its behavior under compositions
and under formation of pointwise sums and products, are thoroughly analyzed in the
first part of \cite{Aichinger-Moosbauer21}.  The next important problem is the
determination of $\mathcal{D}(A,B)$, the set of all functional degrees of maps from
$A$ to $B$. In particular we wish to determine $\delta(A,B)$, the supremum of all
functional degrees of maps from $A$ to $B$.
\\ \indent
In the middle part of their work, in \cite[\S 7-9]{Aichinger-Moosbauer21}, Aichinger
and Moosbauer make progress on this question.  They showed that the
determination of $\delta(A,B)$ for all finite commutative $p$-groups $A$ and $B$ will
give the determination of $\mathcal{D}(A,B)$ for all finite commutative groups $A$
and $B$.    Moreover they recast the determination of $\delta(A,B)$ as an open
problem in commutative algebra \cite[Problem 8.3]{Aichinger-Moosbauer21}.  In the
special case of finite \emph{cyclic} $p$-groups, a solution can then be derived from
work of Wilson \cite{Wilson06}, as Aichinger and Moosbauer point out. As we will see
in \S 4.3, this case is also a consequence of earlier work of Weisman
\cite{Weisman77}.  In order to complete the determination of $\mathcal{D}(A,B)$
when $A$ and $B$ are both finite, it remains to determine $\delta(A,B)$ when $A$
and $B$ are both finite commutative $p$-groups and $A$ is not cyclic.
\\ \\
In this paper, we compute $\delta(A,B)$ for \emph{all} commutative groups $A$ and
$B$: Theorem \ref{MAINTHM1}.  This yields, in particular, an answer to \cite[Problem
8.3]{Aichinger-Moosbauer21}.  We also determine whether $\infty$ lies in
$\mathcal{D}(A,B)$: also in Theorem \ref{MAINTHM1}. As we will see, there are
cases in which $\delta(A,B) = \infty$ but every single map $f: A \ra B$ has finite
functional degree.
\\ \indent
For many pairs $(A,B)$ of commutative groups, it is much easier to see that there is
a function of infinite functional degree than to determine the entire set
$\mathcal{D}(A,B)$.  However, we will compute the full set $\mathcal{D}(A,B)$ for a
class of pairs $(A,B)$ that includes all pairs of finitely generated commutative
groups: Proposition \ref{3.9} and Theorem \ref{3.10}.

\subsection{The Contents of the Paper}
In \S \ref{sec.2} and \ref{sec.3} we give a self-contained exposition of the aspects of
the functional degree that we will need here.  We do so mainly in order to establish
several new intermediate results, some of which will be used later on.\footnote{We
were very taken with the results of \cite{Aichinger-Moosbauer21} and their elegant
exposition of them.  When we can push these results further or contribute a useful
new way to look at them, we do so.  There is more in \cite{Aichinger-Moosbauer21}
than is revisited here, and we will make use of some of their other results in the
sequel \cite{CSII}.}
\\ \indent
Aichinger and Moosbauer place their functional calculus in the setting of maps
between \emph{arbitrary} commutative groups, and throughout this work we try to
follow this lead by working without finiteness restrictions when we can do so.  For
this we need a few tenets of the theory of commutative groups, e.g. the structure of
groups of finite exponent.  This is recalled in \S \ref{sec.2}.  In that section we also
introduce the finite difference calculus and give the definition of functional degree in
those terms.
\\ \indent
Exploiting the fact that the set $B^A$ of all maps from $A$ to $B$ is naturally a
module under the group ring $\Z[A]$, Aichinger and Moosbauer gave an elegant
module-theoretic interpretation of the functional degree. In \S \ref{sec.3} we further
develop the calculus of functional degrees from this module-theoretic point of view.
One unifying feature of our approach is an attention to the effect of the functional
degree of maps $f: A \ra B$ under composition with group homomorphisms
$\varepsilon: A' \ra A$ and $\mu: B \ra B'$: see Lemma \ref{1.3} and Corollary
\ref{1.4}.   We also examine, in \S 3.3,  what happens with the functional degree
when $B$ is a direct product and when $A$ is a direct sum.  In particular we
establish the Diagonalization Theorem (Theorem \ref{DIAGONALTHM}), which
generalizes Lemma 9.3 and Theorem 9.4 of \cite{Aichinger-Moosbauer21}.
\\ \indent
In \S 4 we recall the ideal-theoretic interpretation of $\delta(A,B)$ due to Aichinger
and Moosbauer.  We then establish a key result, Theorem \ref{SUMTHM}. It tells us
how the largest possible degree depends on the direct summands of $A$, if the
\(p\)-group $A$ can be written as a direct sum.  Using this result and the earlier
results, we prove the main results of the paper.

%Second, using their calculus of functional degrees they revisit the work of
%Chevalley and Warning on the number of solutions of a system of polynomial equations
%over a finite field $\F_q$.  If $q = p^N$ then $(\F_q,+) \cong (\Z/p\Z)^N$, and they show
%that in place of polynomial maps from $\F_q^n$ to $\F_q$ of specified degrees, one can
%equally well consider maps $(\Z/p\Z)^{nN} \ra (\Z/p\Z)^N$ of specified functional
%degrees and prove generalizations of Chevalley's Theorem \cite[Thm.
%11.2]{Aichinger-Moosbauer21}, Warning's First Theorem \cite[Thm.
%12.2]{Aichinger-Moosbauer21} and Warning's Second Theorem \cite[Thm.
%14.3]{Aichinger-Moosbauer21}.  It is remarkable that ``ignoring the ring-theoretic
%structure'' in this way leads to deeper results.  In fact their approach can allow
%\emph{more} ring-theoretic generality:
%%Here is such a generalization of the Chevalley-Warning Theorem:
%
%\begin{thm}(Aichinger-Moosbauer \cite[Thm. 12.6]{Aichinger-Moosbauer21})
%\label{INTRO.0}
%Let $R$ be a finite rng\footnote{A \text{rng} is like a ring but not necessarily having a multiplicative identity.  It need not be
%commutative.} of order a power of a prime number $p$. Let $n \in \Z^+$, and let $f_1,\ldots,f_r$ be polynomial
%expressions over $R$ in $N$ variables.  If $\sum_{i=1}^r \deg(f_i) < n$, then
%\[ p \mid \# \{(x_1,\ldots,x_n) \in R^n \mid f_1(x_1,\ldots,x_n) = \ldots = f_r(x_1,\ldots,x_n) = 0 \}. \]
%\end{thm}
%\noindent
%The Chevalley-Warning Theorem is the case $R = \F_q$ of Theorem \ref{INTRO.0}.

\subsection{In the Sequel}

Aichinger and Moosbauer also gave some striking Diophantine applications of their
functional calculus.  Namely, they obtained group-theoretic generalizations
\cite[Thm. 11.1, Thm. 12.2]{Aichinger-Moosbauer21} of the theorems of Chevalley
\cite{Chevalley35} and Warning \cite{Warning35} on systems of polynomials over a
finite field $\F_{p^a}$: the latter result says that for a system of polynomials
$P_1,\ldots,P_r \in \F_{p^a}[t_1,\ldots,t_n]$ with $\sum_{j=1}^r \deg(P_j) < n$, the
number $\# Z$ of simultaneous solutions is divisible by the prime $p$.  In the sequel
\cite{CSII} we will use results from the present paper and a generalization of the
main result of \cite{Wilson06} to give group- and ring-theoretic generalizations of the
results of Ax \cite{Ax64}, Katz \cite{Katz71} and Moreno-Moreno
\cite{Moreno-Moreno95} on higher $p$-adic congruences for $\# Z$.

\section{Finite Differences and the Functional Degree}\label{sec.2}

\subsection{Notation} We denote by $\mathcal{P}$ the set of (positive) prime numbers.
We denote by $\N$ the non-negative integers and put $\Z^+ \coloneqq \N \setminus
\{0\}$.  For $n\in\N$, we set $Z_n:=\Z/n\Z$. In particular, $Z_0=\Z$.  Moreover, we
endow the set
\[ \NN \coloneqq \N \cup \{- \infty,\infty\} \] with the most evident total ordering, in which $-\infty$ is the least element and $\infty$ is the greatest element.
%\\ \\
%For $d \in \Z^+$, we put ${t \choose d} \coloneqq \frac{t(t-1)\cdots(t-d+1)}{d!} \in
%\Q[t]$. For $x \in \N$, we have that ${x \choose d}$ is the usual binomial coefficient
%and is thus a non-negative integer.  Moreover we have ${x \choose d} \in \Z$ for all
%$x \in \Z$: see e.g. \cite[p. 19]{Cahen-Chabert}.   These \textbf{integer-valued
%polynomials} are discussed in \S 4.3 (and in fact a proof that ${t \choose d}$ is
%integer-valued follows from the proof of Theorem \ref{4.3}).  We take ${x \choose 0}:
%\Z \ra \Z$ to be the constant function $1$ and for any negative integer $n$, we take
%${x \choose n}: \Z \ra \Z$ to be the zero function.

\subsection{Preliminaries on Commutative Groups}
Let $(A,+)$ be a commutative group.  For $n \in \Z^+$ we define
\[ A[n] \coloneqq \{x \in A \mid nx = 0\} \]
and
\[ A[n^{\infty}] \coloneqq \{x \in A \mid n^k x = 0 \text{ for some } k \in \Z^+\}. \]
For a prime number $p$, we say that $A$ is a \textbf{p-group} if $A = A[p^{\infty}]$.
\smallskip\\
The \textbf{torsion subgroup} of $A$ is
\[ A[\tors] \coloneqq \{x \in A \mid nx = 0 \text{ for some } n \in \Z^+\} = \bigcup_{n \in \Z^+} A[n]. \]
We say that $A$ is \textbf{torsion} if $A = A[\tors]$.  We say that $A$ has
\textbf{finite exponent} if $A = A[n]$ for some $n \in \Z^+$, in which case the least
such $n$ is called the \textbf{exponent} of $A$ and is denoted by $\exp(A)$.  In the
other case, we may write $\exp(A)=\infty$. We also set
\[ e(A) \coloneqq \begin{cases} \exp(A) & \text{if $\exp(A)<\infty$,} \\ 0 & \text{if $\exp(A)=\infty$.} \end{cases} \]
For a finitely generated commutative group $A$, we denote by $\rank(A)$ the least
$n$ such that $A$ is a direct sum of $n$ cyclic groups.
\smallskip\\
A commutative group $A$ is \textbf{torsion-split} if $A[\tors]$ is a direct summand of
$A$.  Every finitely generated commutative group is torsion-split.   More generally,
Baer showed \cite{Baer36} that if $A[\tors]$ is the direct sum of a group of finite
exponent and a divisible group (i.e., one in which the map $x \mapsto nx$ is
surjective for all $n \in \Z^+$), then $A$ is torsion-split.

\begin{thm}[Pr\"ufer-Baer]
\label{PRUFERBAERTHM} Let $G$ be a commutative group of finite exponent $N$.
Then there is a family $(N_\gamma)_{\gamma\in\Gamma}$ of numbers
$N_\gamma\in\Z^+$ with $N_\gamma \mid N$ for all $\gamma$ and
$\smash{\max\limits_{\gamma\in\Gamma} N_\gamma = N}$, such that $G$ is
isomorphic to the direct sum $\bigoplus\limits_{\gamma \in \Gamma}
Z_{N_\gamma}$.
\end{thm}
\begin{proof}
See e.g. \cite[Thm. 6]{Kaplansky}.
\end{proof}

\subsection{The Functional Degree}
For commutative groups $A$ and $B$, let $B^A$ be the set of all functions $f: A \ra
B$. The set $B^A$ is a commutative group under pointwise addition.  For $f \in B^A$
and $a \in A$, we define the shift operator $\tau_a \in \End B^A$ by
\[ (\tau_a f)(x) \coloneqq f(x+a), \]
and the difference operator $\Delta_a \in \End B^A$ by
\[ (\Delta_a f) \coloneqq \tau_a f - f. \]
In other words, for all $f: A \ra B$ and all $x \in A$ we have
\[ (\Delta_a f)(x) = f(x+a)-f(x). \]
The following is a well-known formula from the calculus of finite differences carried
over to the present context.  The proof is straightforward using induction on $n$.

\begin{lemma}
\label{0.0} Let $a \in A$, $n \in \N$ and let $\Delta_a^n$ be the $n$-fold product
$\Delta_a \cdots \Delta_a \in \End B^A$\!. For all $f \in B^A$ and all $x \in A$,
\[ (\Delta_a^n f)(x) = \sum_{i=0}^n (-1)^i \binom{n}{i} f(x+(n-i)a)
= \sum_{j=0}^n (-1)^{n-j} \binom{n}{j} f(x+ja). \]
\end{lemma}

\begin{dfn}
\noindent For $f \in B^A$ we define the \textbf{functional degree} $\fdeg(f) \in \NN$
as follows:
\begin{itemize}
\item We put $\fdeg(f) = -\infty$ if and only if $f = 0$.\smallskip
\item For $n \in \N$, we say that a nonzero function $f$ has $\fdeg(f) \leq n$ if and
    only if $\Delta_{a_1} \cdots \Delta_{a_{n+1}} f = 0$ for all $a_1,\ldots,a_{n+1} \in
A$.\smallskip
\item If there is an $n \in \N$ such that $\fdeg(f) \leq n$ and $f\neq0$, then we
    define the functional degree $\fdeg(f)$ as the least $n \in \N$ such that $\fdeg(f)
\leq n$.\smallskip
\item If there is no $n \in \N$ such that $\fdeg(f) \leq n$, then we put $\fdeg(f) =
    \infty$.
\end{itemize}
\end{dfn}
\noindent Our definition of functional degree differs from the one given in \cite[\S
2]{Aichinger-Moosbauer21} precisely in that we set $\fdeg(0) = -\infty$, whereas
Aichinger-Moosbauer take $\fdeg(0) = 0$. Our choice is motivated by the
corresponding convention for the degree of polynomials, which ensures that the identity
$\deg(f\cdot g)=\deg(f)+\deg(g)$ for polynomials $f$ and $g$ over an integral
domain holds even when $f=0$ or $g=0$.

\begin{remark}
\label{EX1} Let $f \in B^A$.
\begin{itemize}
\item[a)] (cf. \cite[Lemma 3.1(4)]{Aichinger-Moosbauer21}) We have $\fdeg(f) \leq
    0$ if and only if  $f(x+a) = f(x)$ for all $a,x \in A$, that is, if and only if $f$ is
    constant.\smallskip
\item[b)] \cite[Lemma 3.1(5)]{Aichinger-Moosbauer21} If $f$ is a nonzero group
    homomorphism, then for all $a \in A$ the function $\Delta_a f$ is constant, as for
    all $x \in A$ \[ (\Delta_a f)(x) = f(x+a)-f(x) = f(a) = \text{const.}\] Hence, $\fdeg(f)
    \leq 1$, and since $f$ is nonconstant, $\fdeg(f) = 1$.\smallskip
\item[c)] Suppose $f \in B^A$ has functional degree $1$.  Then for all $a_1,a_2 \in
    A$ we have
\[ 0 = \Delta_{a_1} \Delta_{a_2} f(x) = f(a_1+a_2+x)-f(a_2+x)-f(a_1+x) + f(x). \]
Taking $x = 0$ we get
\[ f(a_1+a_2)-f(0) = (f(a_1)-f(0)) + (f(a_2)-f(0)). \]
So, $x \mapsto f(x)-f(0)$ is a nonzero group homomorphism.
\end{itemize}
\end{remark}

\begin{lemma}
\label{0.1} Let $f \in B^A$.
%\begin{itemize}
%\item[a)]
If $\fdeg(f) = n \in \Z^+$, then for every $a \in A$ we have $\fdeg(\Delta_a f) \leq n-1$
and for some $a \in A$ we have $\fdeg(\Delta_a f) = n-1$.
%\item[b)] If $\Hom(A,B) = 0$, then every nonconstant $f \in B^A$ has infinite functional degree.
%\end{itemize}
\end{lemma}
\begin{proof}
Since $\fdeg(f) \leq n$, for all $a,a_1,\ldots,a_n \in A$ we have
\[ \Delta_{a_1} \cdots \Delta_{a_n} (\Delta_a f) = \Delta_{a_1} \cdots \Delta_{a_n} \Delta_a f = 0, \]
so $\fdeg(\Delta_a f) \leq n-1$.  Since $\fdeg(f)$ is \emph{not} less than or equal to
$n-1$, there are $a_1,\ldots,a_n \in A$ such that
\[ \Delta_{a_1} \cdots \Delta_{a_{n-1}} (\Delta_{a_n} f) = \Delta_{a_1} \cdots \Delta_{a_n} f \neq 0, \]
so $\fdeg(\Delta_{a_n} f)$ is \emph{not} less than or equal to $n-2$.  Thus
$\fdeg(\Delta_{a_n} f) = n-1$.
%If $\fdeg(f) = n$ then $f \in B^A[I^{n+1}] \setminus B^A[I^{n}]$.  For all $a \in A$ we have $\Delta_a \in I$
%thus $\Delta_a(f) \in IB^A[I^{n+1}] \subset B^A[I^n]$, so $\fdeg(\Delta_a f) \leq n-1$.  Since $I^{n}$ is generated by
%$(n)$-fold products of the operators $\Delta_a$, there are $a_1,\ldots,a_n$ such that
%\[ 0 \neq \Delta_{a_1} \cdots \Delta_{a_n} f = \Delta_{a_1} (\Delta_{a_2} \cdots \Delta_{a_{n}} f), \]
%so $\Delta_{a_1} f \notin B^A[I^{n-1}]$ and thus $\fdeg(\Delta_{a_1} f) \geq n-1$.
\end{proof}
\noindent For commutative groups $A$ and $B$, let
\[ \mathcal{D}(A,B) \coloneqq \{ \fdeg(f) \mid f \in B^A\} \subset \NN\]
be the set of all functional degrees of maps between $A$ and $B$.  The main result
of \S 3 is the complete determination of $\mathcal{D}(A,B)$ for certain groups $A$
and $B$, including all finitely generated commutative groups. %\footnote{In what follows, the
%case of most interest to us is the one in which $A$ and $B$ are both finite $p$-groups.}
The set $\mathcal{D}(A,B)$ is completely determined by two pieces of
information: the invariant
\[ \delta^{\circ}(A,B) \coloneqq \sup \{ \fdeg(f) \mid f \in B^A, \ \fdeg(f) < \infty \}, \]
and whether there is a function $f \in B^A$ with $\fdeg(f) = \infty$.  Indeed:

\begin{lemma}
\label{0.1.5} Let $A$ and $B$ be commutative groups.
\begin{itemize}
\item[a)] If $B$ is trivial, then $B^A = \{0\}$, so $\mathcal{D}(A,B) = \{-\infty\}$.
\item[b)] If $A$ is trivial and $B$ is nontrivial, then $\mathcal{D}(A,B) = \{-\infty,0\}$.
\item[c)] If $A$ and $B$ are both nontrivial, then
% and that there is no $f \in B^A$ with $\fdeg(f) =\infty$.  Then
\[\mathcal{D}(A,B) \setminus \{ \infty\} = \{n \in {-\infty} \cup \N  \mid n \leq \delta^{\circ}(A,B)\}. \]
%\item[d)] Supose that $A$ and $B$ are nontrivial and that there is $f \in B^A$ with $\fdeg(f) = \infty$.  Then
%\[\mathcal{D}(A,B) = \{n \in \NN \mid n \leq \delta^{\circ}(A,B)\} \cup \{ \infty\}. \]
\end{itemize}
\end{lemma}
\begin{proof}
a) This is trivial.\smallskip\\
b) In this case, every function is constant (as its domain has only one element), and
some are nonzero.\smallskip\\
c) This follows from Lemma \ref{0.1}.
\end{proof}
\noindent The following consequence of Remark \ref{EX1} and Lemma \ref{0.1.5}
shows a connection between functional degrees and the existence of nontrivial
group homomorphisms from $A$ to $B$.

\begin{cor}
\label{0.1.7} For commutative groups $A$ and $B$, the following are equivalent:
\begin{itemize}
\item[(i)] $\mathcal{D}(A,B) \subset \{-\infty,0,\infty\}$.
\item[(ii)] $\Hom(A,B) = \{0\}$.
\item[(iii)] $\delta^{\circ}(A,B) = 0$.
\end{itemize}
\end{cor}
\noindent Next we introduce a quantity that is closely related to $\delta^{\circ}(A,B)$
but easier to compute:
\[ \delta(A,B) \coloneqq \sup \{\fdeg(f) \mid f \in B^A \}. \]
Unlike $\delta^{\circ}(A,B)$, we can compute $\delta(A,B)$ for all commutative
groups $A$ and $B$: this is because, for large classes of groups $A$ and $B$, we
can simply write down an $f \in B^A$ with $\fdeg(f) = \infty$.
\\ \\
If $\delta(A,B) < \infty$, then $\delta^{\circ}(A,B) = \delta(A,B)$. So it follows from
Lemma \ref{0.1.5} that
\begin{equation*}
%\label{DELTAAB1EQ}
 \mathcal{D}(A,B) = \{n \in \NN \mid n \leq \delta(A,B) \}\ \ \text{if}\ \ \delta(A,B) < \infty.
\end{equation*}

\section{Module-theoretic Interpretation of the Functional Degree}\label{sec.3}

\subsection{Preliminaries on Group Rings}
For a group $(A,+)$ and a ring $R$, let $R[A]$ be the corresponding \textbf{group
ring} \cite[\S 5.6]{Clark-CA}.  Its elements are formal linear combinations $\sum_{a
\in A} r_a [a]$ with $r_a \in R$ and $r_a = 0$ for all but finitely many $a \in A$. The
ring homomorphism
\[ \epsilon: R[A] \ra R,\ \sum_{a \in A} r_a [a] \mapsto \sum_{a \in A} r_a. \]
is called \textbf{augmentation map}, and its kernel is the \textbf{augmentation ideal}
$I$.

\begin{remark}[Reminders on Group Rings]
\label{REMARK1} Let $R$ be a nonzero commutative ring.
\begin{itemize}
\item[a)] The $R$-algebra homomorphism $R[\Z] \ra R[t,t^{-1}]$ that sends $[1]$ to
    $t$ is an isomorphism.  Under this map the augmentation ideal maps to the
    principal ideal $\langle t-1 \rangle$.
\item[b)] Let $N \geq 2$.  The $R$-algebra homomorphism $R[Z_N] \ra
    R[t]/\langle t^N -1 \rangle$ that sends $[1]$ to $t + \langle t^N-1 \rangle$ is an
    isomorphism. Under this map the augmentation ideal maps to the principal ideal
    $\langle t-1 + \langle t^N-1 \rangle \rangle = \langle t-1 \rangle/\langle t^N-1
    \rangle$ (which we identify with $\langle t-1 \rangle \subset R[t]$).
\item[c)] For groups $G_1$ and $G_2$, the natural map
\[ R[G_1] \otimes_R R[G_2] \ra R[G_1 \times G_2], \ [g_1] \otimes [g_2] \mapsto [(g_1,g_2)] \]
is an $R$-algebra isomorphism.
\item[d)] Combining parts b) and c), we get that if $A \cong \bigoplus_{i=1}^r Z_{N_i}$, then
\[ R[A] \cong R[t_1,\ldots,t_r]/\langle t_1^{N_1}-1,\ldots,t_r^{N_r}-1 \rangle. \]
The augmentation ideal is mapped to $\langle t_1-1,\ldots,t_r -1 \rangle/ \langle
t_1^{N_1}-1,\ldots,t_r^{N_r}-1 \rangle$ (which we identify with $\langle
t_1-1,\ldots,t_r-1 \rangle$) under this isomorphism.
\end{itemize}
\end{remark}
\noindent We call $\Z[A]$ the \textbf{integral group ring} of $A$.  The commutative group $B^A$
has the structure of a $\Z[A]$-module via
\[\bigl(\sum_{a \in A} n_a [a]\bigr)f \coloneqq \sum_{a \in A} n_a \tau_a(f). \]
Moreover, if $B$ has finite exponent $n$, then $B$, hence also $B^A$, is a
$Z_n$-module, and the $\Z[A]$-module structure on $B^A$ factors through
$Z_n[A]$. We see that $B^A$ is always a $Z_{e(B)}[A]$-module, even if
$\exp(B)=\infty$.
\begin{lemma}
\label{1.1} $B^A$ is a faithful $Z_{e(B)}[A]$-module.
\end{lemma}
\begin{proof}
Let $0\neq r = \sum_{a \in A} n_a [a]\in Z_{e(B)}[A]$ be arbitrary. Choose $a_0\in
A$ such that $n_{a_0}\neq0$, and let $\hat n_{a_0}\in\Z$ be a least nonnegative
representative of $n_{a_0}\in Z_{e(B)}$. As $0<\hat n_{a_0}<\exp(B)$, we have
$B[\hat n_{a_0}]\neq B$ and there exists an element $b\in B$ with $n_{a_0}b:=\hat
n_{a_0}b\neq0$. From this property it follows easily that if $\delta_{0,b}\in B^A$ is
the function that maps $0$ to $b$ and every other element of $A$ to $0$, then $r
\cdot \delta_{0,b}\neq0$, which proves the claimed faithfulness. Indeed, if we
evaluate the function
\[ r \cdot \delta_{0,b}: x \longmapsto \sum_{a \in A} n_a \delta_{0,b}(a+x) \]
at $x=-a_0$, we see that $[r \cdot \delta_{0,b}](-a_0)=n_{a_0} b\neq0$.
\end{proof}

\begin{remark}[Reminders on Nil Ideals]
\label{REMARK2} Let $R$ be a commutative ring, and let $I$ be an ideal of $R$.
\begin{itemize}
\item[a)] We say that $I$ is \textbf{nil} if every $x \in I$ is nilpotent: i.e., there is $n =
    n(x) \in \Z^+$ such that $x^n = 0$.   We say that $I$ is \textbf{nilpotent} if $I^n =
    0$ for some $n \in \Z^+$.  The \textbf{nilpotency index} $\nu(I)$ of an ideal $I$ is
    the least $n \in \N$ such that $I^n = 0$ or $\infty$ if there is no such $n$.  Thus
    $\nu(I) < \infty$ if and only if $I$ is nilpotent.
\item[b)] Let $x_1,\ldots,x_r \in S$ be nilpotent elements; more precisely, let
    $a_1,\ldots,a_r \in \Z^+$ be such that $x_i^{a_i} = 0$ for all $1 \leq i \leq r$.
    Then $\langle x_1,\ldots,x_r \rangle^{a_1 + \ldots + a_r - (r-1)} = 0$: indeed, if $N
    \geq a_1 + \ldots + a_r - (r-1)$ then in any expression $y = x_1^{b_1} \cdots
    x_r^{b_r}$ with $b_1 + \ldots + b_r = N$ we have $b_i \geq a_i$ for some $i$ and
    thus $y = 0$.  It follows from this that an ideal is nil if and only if it is generated by
    nilpotent elements and also that every finitely generated nil ideal is nilpotent.
% (hence in a Noetherian ring, all nil ideals are nilpotent).
\end{itemize}
\end{remark}
\noindent For a nonzero commutative group $A$ and a nonzero commutative ring
$R$, we denote by $\nu(R[A])$ the nilpotency index of the augmentation ideal $I$ in
$R[A]$.

\begin{example}
\label{EX3} Let $p \in \mathcal{P}$ and let $n \in \Z^+\!$.  According to Remark
\ref{REMARK1}, $\nu(Z_p[Z_{p^n}])$ is the nilpotency index of the ideal $\langle
t-1 \rangle$ in the ring $Z_p[t]/\langle t^{p^n}-1 \rangle$.  Since
$(t-1)^{p^n}=t^{p^n}-1$ in $Z_p[t]$, evidently we have
\[ \nu(Z_p[Z_{p^n}]) = p^n. \]
\end{example}

\begin{lemma}
\label{1.2} Let $A$ be a nontrivial commutative group, and let $R$ be a nonzero
commutative ring.
\begin{itemize}
\item[a)] If $A$ has an element of infinite order, then the augmentation ideal of
    $R[A]$ is not nilpotent.\smallskip
\item[b)] The augmentation ideal of $\Z[A]$ is not nilpotent.\smallskip
\item[c)] Let $m,N \geq 2$.  Suppose that there is $p \in \mathcal{P}$ such that $p
    \mid N$ and $p \nmid m$ and that $A$ has an element of order $m$.  Then
    $\nu(Z_N[A]) = \infty$.\smallskip
\item[d)] If $A$ is finite, then $\nu(R[A]) \geq \max(\exp(A),\rank(A))$.
\end{itemize}
\end{lemma}
\begin{proof}
Let $H$ be a subgroup of $A$.  There is a natural injective ring homomorphism
$\iota: R[H] \hookrightarrow R[A]$. Moreover, if $I_H$ (resp.\ $I_A$) is the
augmentation ideal of $R[H]$ (resp.\ of $R[A])$, then for all $n \in \Z^+$ we have
$\iota(I_H^n) \subset I_A^n$, so $\nu(R[H]) \leq \nu(R[A])$.\smallskip\\
a) Let $H$ be the subgroup generated by an element of $A$ of infinite order. So, $H
\cong (\Z,+)$ and thus $\nu(R[A]) \geq \nu(R[\Z])$. It suffices to show $\nu(R[\Z]) =
\infty$.  The group ring $R[\Z]$ is isomorphic to the Laurent polynomial ring
$R[t,t^{-1}]$, and under this isomorphism the augmentation ideal maps to the
principal
ideal $\langle t-1 \rangle$.  The element $t-1$ is not nilpotent in $R[t,t^{-1}]$.\smallskip\\
b) In view of part a) we may assume that $A$ has an element $x$ of finite order $m
\geq 2$.  Taking $H$ to be the subgroup generated by $x$, we have $I_H \subset
I_A$.  Thus, it is enough to show that the augmentation ideal of $\Z[H] \cong
\Z[Z_m]$ is not nilpotent.   We have
\[ \Z[Z_m] \cong \Z[t]/\langle t^m-1 \rangle, \]
and we claim that the ring $\Z[t]/\langle t^m-1 \rangle$ is reduced, i.e., has no
nonzero nilpotent elements, which suffices.  Since $t^m-1$ is monic, we have an
injective ring homomorphism
\[ \Z[t]/\langle t^m-1 \rangle \hookrightarrow  \Q[t]/\langle t^m-1 \rangle
\cong \prod_{d \mid m} \Q[t]/\langle \Phi_d \rangle, \] where $\Phi_d$ is the $d^\text{th}$
cyclotomic polynomial.  The ring $\prod_{d \mid m} \Q[t]/\langle \Phi_d \rangle$ is a
product of
fields, hence reduced, hence its subring $\Z[t]/\langle t^m-1 \rangle$ is also reduced.\smallskip\\
c) Arguing as above, it suffices to show that $\nu(Z_N[Z_m]) = \infty$.  Let $p
\in \mathcal{P}$ be such that $p \mid N$ and $p \nmid m$.  Then the surjective ring
homomorphism $Z_N \ra Z_p$ induces a surjective ring homomorphism
\[ q: Z_N[Z_m] \ra Z_p[Z_m], \]
and the map $q$ induces for all $n \in \Z^+$ a surjection between the $n$th powers
of the augmentation ideals, so it suffices to show that $\nu(Z_p[Z_m]) = \infty$.
This follows from a similar argument to that of part b): we have
\[ Z_p[Z_m] \cong Z_p[t]/\langle t^m-1 \rangle, \]
and $t^m-1 \in \F_p[t]$ is separable since $p \nmid m$, so $Z_p[t]/\langle t^m-1
\rangle$
is reduced.\smallskip\\
d) If an $a \in A$ has order $n$ then $([a]-[0])^{n-1} \neq 0$. This shows that
$\nu(R[A]) \geq \exp(A)$. If $A$ has rank $n$ and $x_1,\dotsc,x_n\in A$ are such
that $A \cong \bigoplus_{i=1}^n \langle x_i\rangle$, then $\prod_{i=1}^n ([x_i]-[0])
\neq 0$. This shows that $\nu(R[A]) \geq \rank(A)$.
\end{proof}

\subsection{Group Ring Interpretation of the Functional Degree}
With the described $Z_{e(B)}[A]$-module structure on $B^A$, we may view the
operators $\tau_a = [a]$ and $\Delta_a = [a]-[0]$ as elements of the group ring
$Z_{e(B)}[A]$, and observe that the elements $\Delta_a$ lie in the
augmentation ideal $I$.  Alternatively, we may also view $B^A$ as $\Z[A]$-module
and $I$ as the augmentation ideal of $\Z[A]$.  The latter choice allows us to simultaneously consider different codomains $B$.  For a fixed $A$ all
possible $B^A$ are modules over the same ring $\Z[A]$. Generalizing and unifying
both choices, we may also view $B^A$ as $Z_{ke(B)}[A]$-module, for any
$k\in\N$. No matter which ring we choose, the operators $\Delta_a$ always lie in
the corresponding augmentation ideal $I$, and they generate that ideal:

\begin{lemma}
\label{AUGIDEAL} Let $R$ be a commutative ring and let $A$ be a commutative
group.
\begin{itemize}
\item[a)] The augmentation ideal $I$ of $R[A]$ is generated as an $R$-module by
    $\{ \Delta_a \mid a \in A\}$.\smallskip
\item[b)] If $S \subset A$ generates $A$ as a commutative group, then for each $n
    \in \Z^+\!$, the set of $n$-fold products $\{ \Delta_{s_1} \cdots \Delta_{s_n} \mid
    s_1,\ldots,s_n \in S\}$ generates $I^n$ as an ideal.
\end{itemize}
\end{lemma}
\begin{proof}
a) Assume $x = n_1[a_1]+n_2[a_2]+\dotsb+n_N[a_N]$ lies in $I$. We show that $x$
lies in the $R$-submodule $M$ of $I$ spanned by the elements $\Delta_a$. If $N =
1$ then $n_1 = 0$ by the definition of $I$, so $x = 0\in M$.  If $N \geq 2$, then the
element
\begin{eqnarray*}
(n_1+\dotsb+n_N)[a_1]  &=& x - n_2([a_2]-[a_1]) - n_3([a_3]-[a_1]) - \dotsb - n_N([a_N]-[a_1]) \\
  &=& x - n_2(\Delta_{a_2} - \Delta_{a_1}) - n_3(\Delta_{a_3} - \Delta_{a_1})
             - \dotsb - n_N(\Delta_{a_N} - \Delta_{a_1})%\\
%  &\in&M\ \,\subset\,\ I
\end{eqnarray*}
lies in $I$, since $x$ and the $\Delta_{a_i}$ lie in $I$. Hence,
$n_1+n_2+\dotsb+n_N=0$, as in the case $N = 1$. It follows that
\[ x \,=\, n_2(\Delta_{a_2} - \Delta_{a_1}) + n_3(\Delta_{a_3}- \Delta_{a_1})
             + \dotsb + n_N(\Delta_{a_N} - \Delta_{a_1})\,\in\, M. \]
b) We observe that for all $s_1,s_2 \in S$ we have
\[ \Delta_{-s} = [-s]-[0] = -[-s] \cdot ([s]-[0]) = -[-s] \cdot \Delta_s \in \langle \Delta_s \mid s \in S \rangle, \]
 \[ \Delta_{s_1+s_2} = [s_1+s_2]-[0] = [s_1]([s_2]-[0]) + ([s_1]-[0]) = [s_1] \Delta_{s_2} + \Delta_{s_1}
 \in \langle \Delta_s \mid s \in S \rangle. \]
This shows that $\langle \Delta_s \mid s \in S \rangle = \langle \Delta_a \mid a \in A
\rangle = I$; the result for $I^n$ follows immediately.
\end{proof}
\noindent From this Lemma, with $R:=\Z$ or with $R:= Z_{e(B)}$, we obtain
Aichinger and Moosbauer's module-theoretic interpretation of the functional degree:

\begin{lemma}
\label{2.I} %(cf. \cite[Def. 2.1]{Aichinger-Moosbauer21})
Let $A$ and $B$ be commutative groups, and let $I$ be the augmentation ideal of
$\Z[A]$ or of $Z_{e(B)}[A]$. For each $f \in B^A$ and $n \in \N$, the
following statements are equivalent:
\begin{itemize}
\item[(i)] $I^{n+1}$ kills $f$, i.e., $\theta f=0$ for every $\theta\in I^{n+1}$.\smallskip
\item[(ii)] $\fdeg(f) \leq n$.
\end{itemize}
\end{lemma}

\noindent This equivalence can also be expressed by saying that the elements of
$B^A$ with functional degree at most $n$ form the set
$$
 B^A[I^{n+1}] \coloneqq \{f \in B^A \mid I^{n+1}f = 0\}
 = \{f \in B^A \mid\forall\theta\in I^{n+1},\ \theta f = 0\}. %U
$$
The module-theoretic interpretation of the functional degree was introduced in
\cite[Def. 2.1]{Aichinger-Moosbauer21}.  It is a useful perspective, because it
allows us to apply the language and the tools of commutative algebra to the
calculus of finite differences.  Here is one simple example.

\begin{lemma}
\label{2.13} Let $A$ and $B$ be commutative groups.
\begin{itemize}
\item[a)] For $f,g \in B^A$ we have $\fdeg(f+g) \leq
    \max(\fdeg(f),\fdeg(g))$.\smallskip
\item[b)] For all $n \in \NN$ the set
\[ \mathcal{F}_n(A,B) \coloneqq \{f \in B^A \mid \fdeg(f) \leq n \} \]
is a subgroup of $B^A$, as is
\[ \mathcal{F}(A,B) \coloneqq \bigcup_{n < \infty} \mathcal{F}_n(A,B). \]
\end{itemize}
\end{lemma}
\begin{proof}
a) We may assume without loss of generality that $\fdeg(g) \leq \fdeg(f) =: n \in \N$.
Then $I^{n+1}$ kills both $f$ and $g$, so it kills the $\Z[A]$-submodule of $B^A$
generated by $f$ and $g$, so it kills $f+g$.\smallskip\\
b) This follows immediately.
\end{proof}
\noindent Let $\varepsilon: A' \ra A$ and $\mu: B \ra B'$ be homomorphisms of
commutative groups.   This yields group homomorphisms
\[ \varepsilon^*: B^{A} \ra B^{A'}, \ f \mapsto \varepsilon^* f \coloneqq f \circ \varepsilon\]
and
\[ \mu_*: B^A \ra (B')^A, \ f \mapsto \mu_* f \coloneqq \mu \circ f. \]
If $\varepsilon$ is injective (resp.\ surjective), then $\varepsilon^*$ is surjective
(resp.\ injective), while if $\mu$ is injective (resp.\ surjective), then $\mu_*$ is
injective (resp.\ surjective).

\begin{lemma}[Homomorphic Functoriality I]
\label{1.3} Let $A,A',B,B'$ be commutative groups and let $f \in B^A$.  Let
$\varepsilon: A' \ra A$, $\mu: B \ra B'$ be group homomorphisms.
%\[ \varepsilon^* f \coloneqq f \circ \varepsilon \in B^{A'}, \ \mu_* f \coloneqq \mu \circ f \in (B')^A. \]
Then:
\begin{itemize}
\item[a)] We have $\fdeg \varepsilon^* f \leq \fdeg f$, with equality if $\varepsilon$ is
    surjective.\smallskip
\item[b)] We have $\fdeg \mu_* f \leq \fdeg f$, with equality if $\mu$ is injective.
\end{itemize}
\end{lemma}
\begin{proof}
A homomorphism $\varphi: G \ra H$ of commutative groups induces a
homomorphism of group rings $\Z[\varphi]: \Z[G] \ra \Z[H]$. If $I_G$ (resp.\ $I_H$) is
the augmentation ideal of $\Z[G]$ (resp.\ of $\Z[H]$) then, for all $n \in \N$, we have
\[ \Z[\varphi](I_G^n) \subset I_H^n, \]
with equality if $\varphi$ is surjective. \smallskip\\
a) %Let $n:=\fdeg(f)$ and $\sum_{a'\in A'}n_{a'}[a']\in I_{A'}^{n+1}$. Since
%$\Z[\varepsilon](I_{A'}^{n+1})\subset I_{A}^{n+1}$, we know that
%$\Z[\varepsilon]\bigl(\sum_{a'\in A'}n_{a'}[a']\bigr)$ kills $f$. So, for every $y\in A$, we
%have
%$$\sum_{a'\in A'}n_{a'}f(y+\varepsilon(a'))=\bigl(\!\sum_{a'\in A'}n_{a'}[\varepsilon(a')]f\bigr)(y)
%=\bigl(\Z[\varepsilon]\bigl(\!\sum_{a'\in A'}n_{a'}[a']\bigr)f\bigr)(y)=0.
%$$
%Here, we may also choose $y:=\varepsilon(x)$ for any $x\in A'$, and then
%$$\bigl(\bigl(\!\sum_{a'\in A'}n_{a'}[a']\bigr)(\varepsilon^*f)\bigr)(x)
%=\sum_{a'\in A'}n_{a'}f(\varepsilon(x)+\varepsilon(a'))=0.
%$$
%Therefore, $\bigl(\sum_{a'\in A'}n_{a'}[a']\bigr)(\varepsilon^*f)=0$, i.e.\
%$\fdeg(\varepsilon^*f)\leq n=\fdeg(f)$.\smallskip\\
%In the case that $\varepsilon$ is surjective, every $y$ can be written in the form
%$\varepsilon(x)$, and the entire calculation can be turned around to prove that also
%$\fdeg(\varepsilon^*f)\geq \fdeg(f)$.
Using the ring homomorphism $\Z[\varepsilon]$, we may view $B^A$ as a
$\Z[A']$-module, and then $\varepsilon^*: B^A \ra B^{A'}$ is a $\Z[A']$-module
homomorphism. For each $n \in \N$, we have $\Z[\varepsilon](I_{A'}^{n+1}) \subset
I_{A}^{n+1}$ and can conclude as follows: if $I_A^{n+1}$ kills $f \in B^A$, then in
particular $I_{A'}^{n+1}$ kills $f$, hence $I_{A'}^{n+1}$ also kills $\varepsilon^* f$.
This shows that $\fdeg \varepsilon^* f \leq \fdeg f$. In the case that $\varepsilon$ is
surjective, we have that $\varepsilon^*$ is injective, and we can also use that
$\Z[\varepsilon](I_{A'}^{n+1}) = I_{A}^{n+1}$. So, $I_{A'}^{n+1}$ kills $\varepsilon^*
f$ if and only if $I_{A'}^{n+1}$ kills $f$ if and only if $I_A^{n+1}$ kills $f$. Thus
$\fdeg(\varepsilon^*f) = \fdeg(f)$, if $\varepsilon$
is surjective.\smallskip\\
b) The pushforward map $\mu_*: B^A \to (B')^A$ is a $\Z[A]$-module
homomorphism, so if $\fdeg(f) \leq n$ then $I_A^{n+1}$ kills $f$ and thus also kills
$\mu_* f$, so $\fdeg(\mu_* f) \leq n$; we deduce that $\fdeg(\mu_* f) \leq \fdeg(f)$. If
$\mu$ is injective then $\mu_*$ maps $B^A$ onto a $\Z[A]$-submodule of $(B')^A$,
and the annihilator ideal of an element $x$ of a submodule $N$ of a module $M$ is
the same as the annihilator ideal of $x$ viewed as an element of $M$.
\end{proof}
\noindent The following useful result is an immediate consequence of Lemma
\ref{1.3}.

\begin{cor}
\label{1.4} Let $f: A \ra B$ be a map between commutative groups.
\begin{itemize}
\item[a)] (Domain Restriction) Let $\underline{A}$ be a subgroup of $A$, and let
    $f|_{\underline{A}}: \underline{A} \ra B$, $x\mapsto f(x)$ be the \textbf{restriction}
    of $f$ to $\underline{A}$.  Then
\[ \fdeg(f|_{\underline{A}}) \leq \fdeg(f). \]
\item[b)] (Codomain Restriction) Let $\underline{B}$ be a subgroup of $B$ such
    that $f(A) \subset \underline{B}$, and let $f|^{\underline{B}}: A \ra \underline{B}$
    be given by $x\mapsto f(x)$.  Then
\[ \fdeg(f|^{\underline{B}}) = \fdeg(f). \]
    \end{itemize}
\end{cor}
\noindent Lemma \ref{1.3} also implies that if $A,A',B,B'$ are commutative groups
such that $A \cong A'$ and $B \cong B'$ then $\delta(A,B) = \delta(A',B')$.

\subsection{Between Sums and Products}
Next we study maps $f:A \ra B$ when $B$ is a direct product and also when $A$ is a
direct sum.  Since the set-theoretic Cartesian product is also the direct product in the
category of commutative groups but the direct sum in the category of commutative
groups is not the set-theoretic coproduct, one expects the case of maps into a
product to be simpler than the case of maps out of a sum.  This turns out to be true,
but we still have a result for maps out of a sum that is satisfactory for our purposes.

\begin{lemma}[Mappings into Products]
\label{AM3.4} Let $A$ be a commutative group, and let $B = \prod_{\gamma \in
\Gamma} B_\gamma$ be the direct product of commutative groups $B_\gamma$
over a nonempty index set $\Gamma$. For each $\gamma \in \Gamma$, let
$\pi_\gamma: B \ra B_\gamma$ be the canonical projection. For all $f \in B^A$,
we have
\[ \fdeg(f) = \sup_{\gamma \in \Gamma} \fdeg ( \pi_\gamma \circ f). \]
%\begin{itemize}
%
%\item[b)] For $f \in \underline{B}^A$, we have $\fdeg f = \sup_{i \in I} \fdeg (\pi_i \circ f: A \ra B_i)$.
%\end{itemize}
\end{lemma}
\begin{proof}
The set-theoretic identity $B^A = \prod_{\gamma \in \Gamma} B_\gamma^{A}$ is
also a $\Z[A]$-module
isomorphism. %U
For a commutative ring $R$ and a family  $(M_\gamma)_{\gamma \in \Gamma}$ of
$R$-modules, if $f = (f_\gamma) \in \prod_{\gamma \in \Gamma} M_\gamma$ is an
element of the product $R$-module and $J$ is an ideal of $R$, then $J f = 0$ if and
only if $J f_\gamma = 0$ for all $\gamma \in \Gamma$. Thus the least $n \in \N$
such that $I^n$ kills $f$ is the least $n$ such that $I^n$ kills all $f_\gamma$.
%  \\
%b) Since $\underline{B} \subset B$, this follows from part a) and Lemma \ref{AM3.4}b).
\end{proof}
\noindent Now suppose $(A_\gamma)_{\gamma\in \Gamma}$ is a nonempty family
of nontrivial commutative groups $A_\gamma$. We put $A \coloneqq
\bigoplus_{\gamma \in \Gamma} A_\gamma$ and view $A_\gamma$ as a subgroup
of $A$ via the canonical injection $\iota_\gamma: A_\gamma \hookrightarrow A$.
Hence, if $a_\gamma\in A_\gamma$, it makes sense to write $\Delta_{a_\gamma}$
for $\Delta_{\iota_\gamma(a_\gamma)}$.
%Let $R\Z\Z\Z$ be a commutative ring, let $I$ be the augmentation ideal of $R[A]$, and let
%$I_\gamma$ be the augmentation ideal of $R[A_\gamma]$, for each $\gamma\in \Gamma$.
If each $A_\gamma$ is generated by a set $S_\gamma$, we may also write $A =
\bigoplus_{\gamma \in \Gamma} \langle S_\gamma\rangle$.
%$A$ is generated by the set $\bigcup_{\gamma \in \Gamma} S_\gamma$.
With this setup, we can formulate the following result:

\begin{lemma}[Mappings out of Sums]
\label{SUMLEMMA} %Let $\Gamma$ is a nonempty index set.  For all $\gamma \in \Gamma$ let $A_\gamma$ be
%a nontrivial commutative group, and let $S_\gamma$ be any set of generators for $A_\gamma$.
%Put
%\[A \coloneqq \bigoplus_{\gamma \in \Gamma} A_\gamma. \]  For all $\gamma \in \Gamma$, let $\iota_\gamma: A_\gamma \hookrightarrow A$ be the canonical injective homomorphism.  In this way we view $A_\gamma$ as a subgroup of $A$, and for $a_\gamma \in A_\gamma$ we write $\Delta_{a_\gamma}$ for $\Delta_{\iota_\gamma(a_\gamma)}$.
%%With notation as above, let $B$ be a commutative group.
%%Put $R = Z_{e(B)}$.
Assume $A = \bigoplus_{\gamma \in \Gamma} A_\gamma = \bigoplus_{\gamma \in
\Gamma} \langle S_\gamma\rangle$, as described above, and let $B$ be another
commutative group. For each $f \in B^A$ and $d \in \N$, the following
statements are equivalent:
\begin{itemize}
\item[(i)] There are $r\in\N$, elements
    $\gamma_1,\gamma_2,\dotsc,\gamma_r \in \Gamma$ and $d_1,d_2,\dotsc,d_r
    \in \Z^+$ with $d_1 + d_2 + \dotsb + d_r = d$, and finite sequences
%\begin{equation}
%\label{SUMLEMMAEQ1}
\[a_{1,1},\ldots,a_{1,d_1} \in S_{\gamma_1}, \]
\[ \vdots \]
\[a_{r,1},\ldots,a_{r,d_r} \in S_{\gamma_r} \]
%\end{equation}
such that
%\begin{equation}
%\label{SUMLEMMAEQ2}
\[\bigl(\prod_{i=1}^r \prod_{j=1}^{d_i} \Delta_{a_{i,j}}\bigr) f \neq 0. \]
%\end{equation}
\item[(ii)] We have $\fdeg(f) \geq d$.
\end{itemize}
\end{lemma}
\begin{proof}
If $d=0$, the theorem holds, as $\bigl(\prod_{i=1}^r \prod_{j=1}^{d_i}
\Delta_{a_{i,j}}\bigr) f=f$ if $r=0$. So, assume $d\geq1$, and let $I$ be the
augmentation ideal of $\Z[A]$. As $S:=\bigcup_{\gamma \in \Gamma} S_\gamma$ is
a set of generators of $A$, we have that $\{ \Delta_{s_1} \cdots \Delta_{s_d} \mid
s_1,\ldots,s_d \in S\}$ is a set of generators for the ideal $I^d\!$, by Lemma
\ref{AUGIDEAL}. With that, the result follows from the contrapositive of Lemma
\ref{2.I} with $n:=d-1$.
\end{proof}
\noindent
%{\color{red}{This might be a place to make the connection to partial functional degrees.}}
%\\ \\
Our main application of Lemma \ref{SUMLEMMA} will be Theorem \ref{SUMTHM},
which reduces the determination of functional degrees of maps between finite
commutative $p$-groups to the cyclic case.
\\ \\
The next result is a discrete analogue of the fact that a smooth function $f: \R^n \ra
\R^n$ with diagonal Jacobian matrix decomposes as $f = (f_1,\ldots,f_n)$ with $f_i:
\R \ra \R$.   Let
$\Gamma$ be a nonempty index set, for each $\gamma \in \Gamma$ let
$A_\gamma$ and $B_\gamma$ be commutative groups, and put
%\[ A \coloneqq \bigoplus_{\gamma \in \Gamma} A_\gamma
%\subset \prod_{\gamma \in \Gamma} A_\gamma \eqqcolon \tilde{A}
%\ \ \text{and}\ \ B \coloneqq \prod_{\gamma \in \Gamma} B_\gamma.
%\]
%Then we may naturally view each $f_{\bullet} = (f_\gamma)_{\gamma \in \Gamma} \in
%\prod_{\gamma \in \Gamma} B_\gamma^{A_\gamma}$ as a function in
%$B^{\tilde{A}}$, and $\prod_{\gamma \in \Gamma} B_\gamma^{A_\gamma}$ as a
%subgroup of $B^{\tilde{A}}$.  We just need to define
%\[ f_{\bullet}(x) \coloneqq (f_\gamma(x_\gamma))_{\gamma \in \Gamma}
%\quad\text{for all \ $x = (x_\gamma)_{\gamma \in \Gamma} \in \tilde{A}$.} \] These
%$f_{\bullet}$ are exactly the ``diagonal'' functions in $B^{\tilde{A}}$,
\[ A \coloneqq \bigoplus_{\gamma \in \Gamma} A_\gamma
%\subset \prod_{\gamma \in \Gamma} A_\gamma \eqqcolon \tilde{A}
\ \ \ \text{and}\ \ \ B \coloneqq \prod_{\gamma \in \Gamma} B_\gamma.
\]
Then we may naturally view each $f_{\bullet} = (f_\gamma)_{\gamma \in \Gamma} \in
\prod_{\gamma \in \Gamma} B_\gamma^{A_\gamma}$ as a function in $B^A$, and
$\prod_{\gamma \in \Gamma} B_\gamma^{A_\gamma}$ as a subgroup of $B^A$. We
just need to define
\[ f_{\bullet}(x) \coloneqq (f_\gamma(x_\gamma))_{\gamma \in \Gamma}
\quad\text{for all \ $x = (x_\gamma)_{\gamma \in \Gamma} \in A$.} \] These
$f_{\bullet}$ are exactly the ``diagonal'' functions in $B^A$, in the sense that for each
$\gamma \in \Gamma$ the $\gamma$-component of the output depends only on the
$\gamma$-component of the input.

\begin{thm}[Diagonalization Theorem]
\label{DIAGONALTHM} With $A \coloneqq \bigoplus_{\gamma \in \Gamma}
A_\gamma$ and $B \coloneqq \prod_{\gamma \in \Gamma} B_\gamma$, as above,
suppose that for all $\gamma \in \Gamma$ we have
$\Hom\bigl(\bigoplus\limits_{\lambda \neq \gamma} A_\lambda,B_\gamma\bigr) =
\{0\}$. Then
\begin{itemize}
\item[a)] $\fdeg(f_{\bullet}) = \sup\limits_{\gamma \in \Gamma} \fdeg(f_\gamma)$ for
    all $f_{\bullet} = (f_\gamma)_{\gamma \in \Gamma} \in \prod\limits_{\gamma \in
    \Gamma} B_\gamma^{A_\gamma}$.\smallskip
\item[b)] $\mathcal{F}_n(A,B) = \prod\limits_{\gamma \in \Gamma}
    \mathcal{F}_n(A_\gamma,B_\gamma)$ for all $n<\infty$.\smallskip
\item[c)] $\mathcal{F}(A,B) =\bigcup\limits_{n<\infty}\prod\limits_{\gamma \in
    \Gamma} \mathcal{F}_n(A_\gamma,B_\gamma) \subset \prod\limits_{\gamma \in
    \Gamma} \mathcal{F}(A_\gamma,B_\gamma)$.
\end{itemize}
\end{thm}
\begin{proof}
a) For $\gamma \in \Gamma$, let $\iota_\gamma: A_\gamma \hookrightarrow A$ and
$\pi_\gamma: B \ra B_\gamma$ be the canonical maps. Whenever $\lambda \neq
\gamma$ Corollary \ref{0.1.7} implies that $\pi_\gamma \circ f \circ \iota_\lambda$ is
constant, as $\Hom(A_\lambda,B_\gamma) = \{0\}$. It follows (e.g.\ by Lemma
\ref{SUMLEMMA} with $S_\gamma:=A_\gamma$) that
\[\fdeg(\pi_\gamma \circ f_{\bullet}) = \fdeg(\pi_\gamma \circ f_{\bullet} \circ \iota_\gamma)
= \fdeg(f_\gamma),\] for all $\gamma \in \Gamma$. So, by Lemma \ref{AM3.4} we have
\[\fdeg(f_{\bullet}) = \sup_{\gamma \in \Gamma} \fdeg(\pi_\gamma \circ f_{\bullet})
= \sup_{\gamma \in \Gamma} \fdeg(f_\gamma).
\]
b) The inclusion $\mathcal{F}_n(A,B) \supset \prod_{\gamma \in \Gamma}
\mathcal{F}_n(A_\gamma,B_\gamma)$ follows directly from part a). To prove the
opposite inclusion, we first introduce some notation. For each $\gamma \in
\Gamma$, we identify $A=\bigoplus_{\gamma \in \Gamma} A_\gamma$ with
$A_\gamma \times \bigoplus_{\lambda \neq \gamma} A_\lambda$.  Having done so, each $x_\gamma \in A_\gamma$ and  $x^\gamma \in \bigoplus_{\lambda
\neq \gamma} A_\lambda$ yields an element $(x_\gamma,x^\gamma)$ of $A$. Now, let $f
\in \mathcal{F}_n(A,B)$. We need to show that $f \in \prod_{\gamma \in \Gamma}
\mathcal{F}_n(A_\gamma,B_\gamma)$. Consider first just one fixed $\gamma \in
\Gamma$ and one fixed $x_\gamma \in A_\gamma$. Let $f^\gamma:A\ra
B_\gamma$ be the $\gamma$-component of $f$.
%We want to show that for all $x^\gamma,b^\gamma \in \bigoplus_{\lambda \neq \gamma} A_\lambda$
%we have $f_\gamma(x_\gamma,x^\gamma) = f_\gamma(x_\gamma,b^\gamma)$.
As $f$ has finite functional degree, so does the function
\[ g_{\gamma,x_\gamma}: \bigoplus_{\lambda \neq \gamma} A_\lambda \ra B_\gamma\ ,
\quad g_{\gamma,x_\gamma}(x^\gamma) := f^\gamma(x_\gamma,x^\gamma). \]
Since $\Hom(\bigoplus_{\lambda \neq \gamma} A_\lambda,B_\gamma) = \{0\}$,
Corollary \ref{0.1.7} implies that $g_{\gamma,x_\gamma}$ is constant,
%.  This means that for all $x^\gamma,y^\gamma \in \bigoplus_{\lambda \neq
%\gamma} A_\lambda$ we have $f^\gamma(x_\gamma,x^\gamma) =
%f^\gamma(x_\gamma,y^\gamma)$,
and we may write $f_\gamma(x_\gamma)$ for that constant value. This applies to all
$x_\gamma \in A_\gamma$ and all $\gamma \in \Gamma$. Now, because $f$ lies in
$\mathcal{F}_n(A,B)$, the functions $f_\gamma:x_\gamma\mapsto
f_\gamma(x_\gamma)$ lie in $\mathcal{F}_n(A_\gamma,B_\gamma)$, and $f =
(f_\gamma) \in \prod_{\gamma \in \Gamma} \mathcal{F}_n(A_\gamma,B_\gamma)$,
indeed.\smallskip\\
c) This follows from part b), as $\prod_{\gamma \in \Gamma}
    \mathcal{F}_n(A_\gamma,B_\gamma) \subset \prod_{\gamma \in \Gamma}
    \mathcal{F}(A_\gamma,B_\gamma)$ for all $n<\infty$.
%$$\mathcal{F}_n(A,B)=\bigcup_{n<\infty}\mathcal{F}_n(A,B)=\bigcup_{n<\infty}\prod\limits_{\gamma
%\in \Gamma} \mathcal{F}_n(A_\gamma,B_\gamma) \subset \prod\limits_{\gamma \in \Gamma}
%\mathcal{F}(A_\gamma,B_\gamma).$$
\end{proof}

\begin{remark}
If we define $f_{\bullet}(x)$ only for $x\in \bigoplus_{\gamma \in \Gamma}
A_\gamma$, as we did, Theorem \ref{DIAGONALTHM} holds also with the direct
sum $\bigoplus_{\gamma \in \Gamma} B_\gamma$ in the place of the direct product
$B:=\prod_{\gamma \in \Gamma} B_\gamma$. One can also formulate that theorem
for the direct product $\prod_{\gamma \in \Gamma} A_\gamma$ instead of the direct
sum $\bigoplus_{\gamma \in \Gamma} A_\gamma$, but this is not possible in Lemma
\ref{SUMLEMMA}, and the direct sum is also more important in the study of torsion
groups.
\end{remark}

\noindent Our main application of Theorem \ref{DIAGONALTHM} is to the case when
$A$ and $B$ are both torsion group.  We take $\Gamma$ to be the set of
prime numbers $\mathcal{P}$, and decompose $A$ and $B$ as
\[ A = \bigoplus_{p \in \mathcal{P}} A[p^{\infty}]\ \ \text{and}\ \ B = \bigoplus_{p \in \mathcal{P}} B[p^{\infty}]
\subset \prod_{p \in \mathcal{P}} B[p^{\infty}] \eqqcolon \tilde{B}.
\]
We may apply Theorem \ref{DIAGONALTHM} in this setting.
%with $\tilde{B}$ in the place of $B$.
Moreover, if $B[p^{\infty}]=\{0\}$ for all but finitely many $p \in \mathcal{P}$, e.g.\ if
$B$ has finite exponent, then $\tilde{B}=B$ and we can write all the direct products
in the theorem as direct sums. In that situation, we can then also use that
$$\bigcup_{n<\infty}\bigoplus_{p \in \mathcal{P}}
\mathcal{F}_n(A[p^{\infty}],B[p^{\infty}]) = \bigoplus_{p \in \mathcal{P}}
\mathcal{F}(A[p^{\infty}],B[p^{\infty}])$$ and obtain the following corollary:
\begin{cor}
\label{DIAGONALCOR} If $A$ and $B$ are torsion groups, and if $B[p^{\infty}]=\{0\}$
for all but finitely many $p \in \mathcal{P}$, then\smallskip
\begin{itemize}
\item[a)] $\fdeg(f_{\bullet}) = \max\limits_{p \in \mathcal{P}} \fdeg(f_p)$ for all
    $f_{\bullet} = (f_p)_{p \in \mathcal{P}} \in \bigoplus\limits_{p \in \mathcal{P}}
    B[p^{\infty}]^{A[p^{\infty}]}$,\smallskip
\item[b)] $\mathcal{F}_n(A,B) = \bigoplus\limits_{p \in \mathcal{P}}
    \mathcal{F}_n(A[p^{\infty}],B[p^{\infty}])$ for all $n<\infty$,\smallskip
\item[c)] $\mathcal{F}(A,B) = \bigoplus\limits_{p \in \mathcal{P}}
    \mathcal{F}(A[p^{\infty}],B[p^{\infty}])$.
\end{itemize}
\end{cor}
%%\begin{proof}
%%First part b) and then part a) follow directly from the corresponding parts of the
%%theorem (with the new $\tilde{B}$ in the place of the old $B$), by intersecting those
%%results with $B^A$ as subset of $\tilde{B}^A$. We only need the insight that
%%$$B^A\cap\prod_{p \in \mathcal{P}}B[p^{\infty}]^{A[p^{\infty}]}=\bigoplus_{p \in \mathcal{P}}
%%    B[p^{\infty}]^{A[p^{\infty}]},$$
%%$$B^A\cap\prod_{p \in \mathcal{P}}
%%    \mathcal{F}(A[p^{\infty}],B[p^{\infty}])=\bigoplus_{p \in \mathcal{P}}
%%    \mathcal{F}(A[p^{\infty}],B[p^{\infty}])$$ and
%%$$B^A\cap\mathcal{F}(A,\tilde{B})=\mathcal{F}(A,B).$$
%%\end{proof}
\noindent We will use the previous results to reduce the study of maps of finite
functional degree between torsion groups to the $p$-primary case.

\section{The set $\mathcal{D}(A,B)$ of functional degrees}

\subsection{Ideal-theoretic Interpretation of $\delta(A,B)$}
Recall that for commutative groups $A$ and $B$, we have defined $\delta(A,B)$ to
be the supremum of all functional degrees of maps $f \in B^A$ and
$\mathcal{D}(A,B)$ to be the set of all functional degrees of maps $f \in B^A$.  In
this section we will compute $\delta(A,B)$ for all $A$ and $B$ and
$\mathcal{D}(A,B)$ for a class of commutative groups that includes all finitely
generated groups $A$ and $B$.
\\ \\
For $f \in B^A$ and $n \in
\N$, by Lemma \ref{2.I}, we have
\[\fdeg(f) = n \iff  f \in B^A[I^{n+1}] \setminus B^A[I^n]. \]
Considering this connection between the degree of $f$ and the powers of the
augmentation ideal $I$ of $Z_{e(B)}[A]$, it is no surprise that the quantity
$\delta(A,B)$ has an ideal-theoretic interpretation:

\begin{thm}
\label{1.5} Let $A$ and $B$ be nontrivial commutative groups. Then $$\delta(A,B) =
\nu(Z_{e(B)}[A])-1.$$
\end{thm}
\begin{proof}
Let $I$ be the augmentation ideal of $Z_{e(B)}[A]$, and put $\nu \coloneqq
\nu(Z_{e(B)}[A])$.  Lemma \ref{1.1} tells us that $B^A$ is a faithful
$Z_{e(B)}[A]$-module.  If $\nu = \infty$, then for all $n \in \Z^+$ there is $\eta  \in
I^n \setminus \{0\}$. By faithfulness, there is $f \in B^A$ such that $\eta f \neq 0$,
and it follows that $\delta(A,B) = \infty$.  If $\nu < \infty$ then $I^{\nu} = 0$, so $B^A
= B^A[I^{(\nu-1)+1}]$, so $\delta(A,B) \leq \nu-1$.  The converse is the same as
above: there is $\eta \in I^{\nu-1} \setminus\{0\}$, hence by faithfulness there is $f \in
B^A$ such that $I^{\nu-1} f \neq 0$, so $\delta(A,B) \geq \nu-1$.
\end{proof}
\noindent We deduce:
\begin{cor}
The following are equivalent:
\begin{itemize}
\item[(i)] $\delta(A,B) < \infty$.
\item[(ii)] The augmentation ideal of $Z_{e(B)}[A]$ is nilpotent.
\end{itemize}
\end{cor}
\begin{proof}
This is immediate from the theorem and Remark \ref{REMARK2}a).
\end{proof}
\noindent We also see that:
\begin{cor}
\label{1.5COR} For commutative groups $A$ and $B$, the quantity $\delta(A,B)$
depends only on $A$ and on the exponent $\exp(B)$ of $B$.
\end{cor}
\noindent Thus we have an interplay between the ideal theory of the group ring
$Z_{e(B)}[A]$ and the structure of its faithful module $B^A$.  Notice that the
nilpotency index of any ideal $J$ in a commutative ring $R$ can be computed using
any faithful $R$-module $M$: it is the least $n\in \Z^+$ such that $M = M[J^n]$ or
$\infty$ if no such $n$ exists.  Because of this, having a ``concretely given'' faithful
$R$-module can be useful for studying the ideal theory of $R$.  Following
Aichinger-Moosbauer \cite[\S 7]{Aichinger-Moosbauer21} we introduce a natural
family of elements of $B^A$, the \textbf{delta functions}, that in many cases can be
shown (via their generation properties) to be elements of $B^A$ of maximal
functional degree.  They provide a convenient tool for giving lower bounds on
$\delta(A,B)$.

\subsection{Delta Functions}
Let $A$ and $B$ be commutative groups.  For $a \in A$ and $b \in B$, we define the
\textbf{delta function} $\delta_{a,b} \in B^A$ via
\[ \delta_{a,b}(x) \coloneqq
\begin{cases} b & \text{if } x = a, \\0 & \text{if } x \neq a. \end{cases} \]
\begin{prop}
\label{DELTAPROP} Let $A$ and $B$ be nontrivial commutative groups.
\begin{itemize}
\item[a)] The $\Z[A]$-submodule of $B^A$ generated by $\mathcal{D} \coloneqq
    \{\delta_{0,b} \mid b \in B\}$ is the set of all functions $f: A \ra B$ such that $f(a)
    = 0$ for all but finitely many $a \in A$.  Thus the subset $\mathcal{D}$
    generates $B^A$ as a $\Z[A]$-module if and only if $A$ is finite.\smallskip
\item[b)] If $A$ is finite and $B= Z_b$, for some $b \geq 2$, then the
    $Z_b[A]$-module $B^A$ is free of rank $1$ with $\delta_{a,u}$ as generator, for
    any choice of $a \in A$ and $u \in Z_b^{\times}$.\smallskip
\item[c)] If $A$ is finite and $B$ has finite exponent then, for each $a \in A$ and
    each $b \in B$ of order $\exp(B)$,
\[ \delta(A,B) = \fdeg(\delta_{a,b}). \]
%\item[d)] If $A$ is infinite then, for each nonzero $b \in B$,
%    $$\fdeg(\delta_{0,b}) = \infty.$$
\end{itemize}
\end{prop}
\begin{proof}
a) For $a \in A$ we have that $\tau_a \delta_{0,b} = \delta_{a,b}$, from which the result
follows easily.\smallskip \\
b) For $B := Z_b^{\times}$, the $\Z$-module generated by $\delta_{0,u}$
contains $\mathcal{D}$, so by part a) the $\Z[A]$-module generated by
$\delta_{0,u}$ is $Z_b^A$.  Being a faithful $Z_b[A]$-module that is
generated by $\delta_{0,u}$,
the module $Z_b^A$ is therefore free of rank $1$.\smallskip \\
c) Let $U$ be the set of elements of $B$ of order $\exp(B)$.  We claim that $U$
generates $B$ as a group.  To see this, let $\underline{B}$ be the subgroup
generated by $U$. We first suppose that $B = B[p^{\infty}]$ for some $p \in
\mathcal{P}$.  In this case, if $x \in U$ and $y \in B \setminus U$ then $x+y \in U$.
It follows that $y \in \underline{B}$, i.e.\ $\underline{B}=B$. In the general case, for
$u \in U$ we write $u = \sum_{p \in \mathcal{P}} u_p$ with $u_p \in B[p^{\infty}]$.
For each $p \in \mathcal{P}$ the element $u_p$ lies in the cyclic subgroup
generated by $u$, hence also in $\underline{B}$, and moreover $u_p$ has maximal
order $\exp(B[p^{\infty}])$.  Conversely, every element $v$ of $B[p^{\infty}]$ of
maximal order is of the form $u_p$ for some $u \in U$.  This shows that
$\underline{B}$ contains $B[p^{\infty}]$ for all $p \in \mathcal{P}$, so $\underline{B}
= B$,
indeed.\smallskip\\
% indeed,fix $x \in U$.
%Then for every $y \in B \setminus \langle x \rangle$ we have $x+y \in U$, so $y$ lies in the subgroup generated by $U$.
Since $\delta_{0,u_1} + \delta_{0,u_2} = \delta_{0,u_1+u_2}$, it follows that
\[\mathcal{U} \coloneqq \{ \delta_{0,u} \mid u \text{ has order } \exp(B)\} \]
generates $B^A$ as a $Z_{\exp(B)}[A]$-module.
%Next we claim that $\fdeg(\delta_{0,u_1}) = \fdeg(\delta_{0,u_2})$ for all $u_1,u_2 \in U$.  Indeed,
Moreover, all $\delta_{0,u}$ have the same degree. This is because for each $u\in
U$ the functional degree of $\delta_{0,u}: A \ra B$ is equal to the functional degree
of $\delta_{0,u}:A \ra \langle u \rangle$, by Corollary \ref{1.4}, and thus also equal
to the functional degree of $\delta_{0,1}: A \ra Z_{\exp(B)}$, by Lemma \ref{1.3},
which does not depend on $u$. It follows that the least power of $I$ (if any) that
kills $\delta_{0,b}$ is the least power of $I$ that kills $\mathcal{U}$, which is the
least power of $I$ that kills $B^A$. Hence, $\delta(A,B) = \fdeg(\delta_{0,b}) =
\fdeg(\delta_{a,b})$.%\smallskip \\
%d) \emph{Case 1, $A$ contains an element $a_0$ of infinite order:} In this case, for
%arbitrarily large $n\in\Z^+$, the elements $0a_0,1a_0,\dotsc,na_0$ are pairwise
%different, so that
%\[ [\Delta_{a_0}^n\delta_{0,b}](0)
%= \sum_{j=0}^n (-1)^{n-j} {n \choose j} \delta_{0,b}(0+ja_0)
%=(-1)^n\delta_{0,b}(0)=(-1)^nb\neq0, \] as in Lemma \ref{0.0}. So, $\fdeg(\delta_{a,b})
%= \fdeg(\delta_{0,b})=\infty$.\smallskip\\
%\emph{Case 2, $A$ contains no element of infinite order:} In this case, for arbitrarily
%large $n\in\Z^+$, one can find elements $a_1, a_2, \dotsc, a_n\in A$ such that
%$a_{j+1}\not\in\langle a_1, a_2, \dotsc, a_j\rangle$, for $j=1,2,\dotsc,n-1$. This
%yields $2^n$ pairwise different elements of the form $a_I:=\sum_{i\in I}a_i$ with
%$I\subset \{1,2,\dotsc,n\}$. It follows that
%$$[\Delta_{a_n}\Delta_{a_{n-1}}\dotsb\Delta_{a_1}\delta_{0,b}](0)
%=\sum_{I\subset \{1,\dotsc,n\}}\!\!(-1)^{n-|I|}\delta_{a_I,b}(0) =(-1)^n\delta_{0,b}(0)
%=(-1)^nb\neq0.
%$$
%So, $\fdeg(\delta_{a,b}) = \fdeg(\delta_{0,b})=\infty$, as before.
\end{proof}

\subsection{Computing $\delta(Z_{p^{\alpha}},Z_{p^{\beta}})$}
We will now compute  $\delta(Z_{p^{\alpha}},Z_{p^{\beta}})$ by two different
arguments. The first uses delta functions, via a 1977 result of Weisman.  The
second uses group rings, via a 2006 result of Wilson.
\\ \\
For $n \in \Z^+$ and $j,k \in \Z$, we put
\[ M_k(j,n) \coloneqq \sum_{\atop{0 \leq i \leq n}{i \equiv j \Mod{k}\!\!\!\!\!\!\!\!\!\!}} (-1)^i \binom{n}{i}. \]

\begin{thm}[Weisman]
\label{WEISMANTHM} Let $p \in \mathcal{P}$, $\alpha,\beta \in \Z^+$ and $j \in
\Z$.
\begin{itemize}
\item[a)] If $n \geq \beta((p-1)+1)p^{\alpha-1}$, then $M_{p^\alpha}(j,n) \equiv 0
    \pmod{p^\beta}$.\smallskip
\item[b)] $M_{p^\alpha}(j,\beta((p-1)+1)p^{\alpha-1}-1) \equiv (-p)^{\beta-1}
    \pmod{p^\beta}$.
\end{itemize}
\end{thm}
\begin{proof}
This is the main result of \cite{Weisman77}.
\end{proof}

\begin{lemma}[Wilson]
\label{Wilson} Let $p \in \mathcal{P}$ and let $\alpha,\beta \in \Z^+$\!.  Then
\begin{equation*}
%\label{WILSONEQ1}
 (t-1)^{(\beta(p-1)+1)p^{\alpha-1}-1} \equiv (-p)^{\beta-1} \sum_{i=0}^{p^{\alpha}-1} t^i
 \pmod{t^{p^\alpha}\!\!-1, p^\beta }.
\end{equation*}
\end{lemma}
\begin{proof}
This is \cite[(22)]{Wilson06}\footnote{In \cite{Wilson06} Wilson speaks of a ``sketch
of a derivation,'' but in fact he provides a complete proof.}.
\end{proof}

\begin{thm}
\label{WEISMANWILSON} Let $p \in \mathcal{P}$, let $\alpha,\beta \in \Z^+\!$ and
let $B$ be a commutative group of exponent $p^\beta$.  Then
\begin{equation*}
%\label{WILSONEQ2}
\delta(Z_{p^{\alpha}},B) = (\beta(p-1)+1)p^{\alpha-1}-1.
\end{equation*}
\end{thm}
\begin{proof}
(Via Weisman) By Corollary \ref{1.5COR}, we may assume that $B=Z_{p^{\beta}}$.
By Proposition \ref{DELTAPROP}b), then
$\delta(Z_{p^{\alpha}},B)=\fdeg(\delta_{0,1})$.  So, by Lemma \ref{AUGIDEAL} and
Lemma \ref{2.I}, we need to determine the least $n \in \N$ such that $\Delta^n
\delta_{0,1}=0$, where $\Delta:=\Delta_1$.  If we view $\delta_{0,1}$ as map from
$Z_{p^{\alpha}}$ into $\Z$, this condition is met if and only if
\[\forall x \in Z_{p^{\alpha}}, \ (\Delta^n \delta_{0,1})(x) \equiv 0 \pmod{p^\beta}.\]
However, Lemma \ref{0.0} tells us that
\[ (\Delta^n \delta_{0,1})(x) =
\sum_{\atop{0 \leq j \leq n}{\!\!\!j \equiv -x \Mod{p^\alpha}\!\!\!\!\!\!\!\!\!\!\!\!\!}} (-1)^{n-j}
\binom{n}{j} = \pm M_{p^\alpha}(-x,n), \] and by Theorem \ref{WEISMANTHM} we know
exactly when this is zero modulo $p^\beta$. We see that
\[ \delta(Z_{p^{\alpha}},B) = (\beta(p-1)+1)p^{\alpha-1}-1. \qedhere \]
\end{proof}

\begin{proof}
(Via Wilson) Theorem \ref{1.5} gives \[\delta(Z_{p^{\alpha}},B) =
\nu(Z_{p^{\beta}})[Z_{p^{\alpha}}])-1. \] Therefore, by Remark \ref{REMARK1}b), it
suffices to show that the nilpotency index of the augmentation ideal $\langle t-1
\rangle$ in the ring $Z_{p^{\beta}}[t]/\langle t^{p^\alpha}-1 \rangle$ is
$(\beta(p-1)+1)p^{\alpha-1}$. Phrasing this in terms of the ring $\Z[t]$, we wish to
show that the least $N \in \Z^+$ such that $(t-1)^N$ lies in the ideal $J \coloneqq
\langle p^\beta,t^{p^{\alpha}}-1 \rangle$ is $(\beta(p-1)+1)p^{\alpha-1}$.  This
follows from Lemma \ref{Wilson}. Indeed, that congruence directly entails
$(t-1)^{(\beta(p-1)+1)p^{\alpha-1}-1}
\notin J$ , and indirectly (after multiplying both sides by $t-1$%
%and summing the finite geometric series
) implies $(t-1)^{(\beta(p-1)+1)p^{\alpha-1}} \equiv 0 \pmod{t^{p^\alpha}\!\!-1,
p^\beta }$ and thus $(t-1)^{(\beta(p-1)+1)p^{\alpha-1}} \in J$.
\end{proof}

\subsection{The $p$-Primary Sum Theorem}

\begin{thm}[$p$-Primary Sum Theorem]
\label{SUMTHM} Let $p \in \mathcal{P}$.  For $1 \leq i \leq r$, let $A_i$ be a
nonzero finite commutative $p$-group, let $A \coloneqq \bigoplus_{i=1}^r A_i$, and
let $B$ be a commutative group of exponent $p^\beta$.  Then
\[ \delta \left( A,B\right) = \max_{\underline{\beta}}\, \sum_{i=1}^r \delta(A_i,Z_{p^{\beta_i+1}}), \]
where the maximum extends over all $\underline{\beta} = (\beta_1,\ldots,\beta_r) \in
\N^r$ with $\beta_1 + \ldots + \beta_r = \beta-1$.
\end{thm}
\begin{proof}
%We prove this using a sequence of equivalences.
By Corollary \ref{1.5COR} we may assume that $B = Z_{p^{\beta}}$. To be able to
apply Proposition \ref{DELTAPROP}c), we further define for $1 \leq i \leq r$ and
$\ell\in\Z^+$ the delta functions $$\delta^i_\ell\coloneqq\delta_{0,1} \in
Z_{p^{\ell}}^{\,A_i},\ \ \delta^i\coloneqq\delta_{0,1} \in \Z^{A_i}\ \ \text{and}\ \
\delta_\ell \coloneqq \delta_{0,1} \in Z_{p^{\ell}}^{\,A}.
$$
To use Lemma \ref{SUMLEMMA} (with the $A_i$ in the place of the generating sets $S_i$), we write
$\underline{a}\in\prod_{i=1}^r A_i^{d_i}$ to say that $\underline{a}$ is a family
$(a_{i,j})$ with $d_i$ entries $a_{i,1},a_{i,2},\dotsc,a_{i,d_i}$ in $A_i$, for each
$1\leq i\leq r$. Similarly, $\underline{x}\in\prod_{i=1}^r A_i$ means that
$\underline{x}$ is an $r$-tuple with $i^{\text{th}}$ entry $x_i$ in $A_i$, etc. Using the
map $\sigma_r:\N^r\to\N$, $\underline{d}\mapsto\sum_{i=1}^rd_i$ we will also write
$\sigma_r^{-1}(d)$ to denote the set of all $\underline{d}\in\N^r$ whose sum of
entries is $d$, for a given $d\in\N$. Moreover, we use the $p$-adic valuation $v_p:
\Z \ra \N \cup \{\infty\}$ with $v_p(0) = \infty$ and $v_p(n) = \max \{a\in\Z^+\mid p^a\
\text{\small divides}\ n\}$ for $n \in \Z \setminus \{0\}$. Now let $d\in\N$. It is enough
to prove the equivalence
$$\delta(A,Z_{p^{\beta}})\geq d\ \ \Longleftrightarrow\ \ \max_{\underline{\beta}}
\sum_{i=1}^r \delta(A_i,Z_{p^{\beta_i+1}})\geq d.
$$
We show this through the following chain of equivalences (where the range of the
variables is specified the first time they appear but not thereafter):
\begin{align*}
   & & &\delta(A,Z_{p^{\beta}})\geq d \\[4pt]
\Leftrightarrow & & &\fdeg(\delta_\beta)\geq d \\[1pt]
\Leftrightarrow & &
\exists\underline{d}\in\sigma_r^{-1}(d):\,\exists\underline{a}\in\prod_{i=1}^r A_i^{d_i}:\,\exists x\in A:\ \
&  \Bigl( \prod_{i=1}^r \prod_{j=1}^{d_i} \Delta_{a_{i,j}} \delta_\beta\Bigr)(x)\neq 0 \in Z_{p^{\beta}} \\
\Leftrightarrow & &
\exists\underline{d}:\,\exists\underline{a}:\,\exists\underline{x}\in\prod_{i=1}^r A_i:\ \
&  \prod_{i=1}^r \Bigl( \bigl( \prod_{j=1}^{d_i} \Delta_{a_{i,j}}\delta^i_\beta\bigr)(x_i)\Bigr)\neq 0 \in Z_{p^{\beta}} \\
\Leftrightarrow & &
\exists\underline{d}:\,\exists\underline{a}:\,\exists\underline{x}:\ \
&  \sum_{i=1}^r v_p\Bigl(\bigl(\prod_{j=1}^{d_i}\Delta_{a_{i,j}}\delta^i\bigr)(x_i)\Bigr)\leq \beta-1 \\
\Leftrightarrow & &
\exists\underline{d}:\,\exists\underline{a}:\,\exists\underline{x}:\,
\exists\underline{\beta}\in\sigma_r^{-1}(\beta-1):\,\forall 1\leq i\leq r:\ \ %\forall i\in\{1,\dotsc,r\}:\ \
&  v_p\Bigl(\bigl(\prod_{j=1}^{d_i}\Delta_{a_{i,j}}\delta^i\bigr)(x_i)\Bigr)\leq \beta_i \\
\Leftrightarrow & &
\exists\underline{d}:\,\exists\underline{\beta}:\,\forall i:\,
 \exists\underline{a_i}\in A_i^{d_i}:\,\exists x_i\in A_i:\ \
&  \bigl(\prod_{j=1}^{d_i}\Delta_{a_{i,j}}\delta^i_{\beta_i+1}\bigr)(x_i) \neq 0 \in Z_{p^{\beta_i+1}} \\[1pt]
\Leftrightarrow & &
\exists\underline{\beta}:\,\exists\underline{d}:\,\forall i:\ \
&  \fdeg(\delta^i_{\beta_i+1})\geq d_i \\[2pt]
\Leftrightarrow & &
\exists\underline{\beta}:\ \
&  \sum_{i=1}^r \fdeg(\delta^i_{\beta_i+1})\geq d \\
\Leftrightarrow & &
&  \max_{\underline{\beta}}\sum_{i=1}^r \delta(A_i,Z_{p^{\beta_i+1}})\geq d. \qedhere
\end{align*}
\end{proof}

\subsection{Computation of $\delta(A,B)$}
The following result computes $\delta(A,B)$ for all nontrivial commutative groups $A$ and $B$,
answering a question of Aichinger-Moosbauer. When $\delta(A,B) = \infty$ we
specify whether every $f \in B^A$ has finite functional degree or whether there are
functions of degree $\infty$.

\begin{thm}
\label{MAINTHM1} \label{3.8} Let $A$ and $B$ be nontrivial commutative groups.
\begin{itemize}
\item[a)] If $A$ is infinite, then $\fdeg(\delta_{0,b}) = \infty$ for all $b \in B
    \setminus \{0\}$, in particular, $\delta(A,B) = \infty$.\smallskip
\item[b)] If there is no $p \in \mathcal{P}$ such that $A$ is a finite $p$-group and
    $B$ is a $p$-group, then $\fdeg(\delta_{0,b}) = \infty$ for some $b \in B$, in
    particular, $\delta(A,B) = \infty$.\smallskip
\item[c)] If, for some $p \in \mathcal{P}$, $A$ is a finite $p$-group, say $A \cong
    \bigoplus_{i=1}^r Z_{p^{\alpha_i}}$ with $\alpha_1 \geq \ldots \geq \alpha_r$,
    and $B$ is a $p$-group of finite exponent $p^\beta$, then
\begin{equation*}
\delta(A,B) = \sum_{j=1}^r (p^{\alpha_j}-1) + (\beta-1)(p-1)p^{\alpha_1-1}.
\end{equation*}
\item[d)] If, for some $p \in \mathcal{P}$, $A$ is a finite $p$-group and $B$ is a
    $p$-group of infinite exponent, then every $f \in B^A$ has finite functional degree
    but $\delta(A,B) = \infty$.
\end{itemize}
\end{thm}
\begin{proof}
a) Let $b \in B\setminus \{0\}$. Since $\delta_{0,b}(A)=\{0,b\} \subset \langle b
\rangle$, Corollary \ref{1.4}b) reduces us to the case $B = \langle b \rangle$. We
may even assume $B= Z_p$ and $b=1\in Z_p$, for a suitable $p\in\mathcal{P}$.
This is because there is always a prime $p$ such that a homomorphism $\mu:
\langle b \rangle \ra Z_q$ with $\mu(b)=1$ exists, and then Lemma \ref{1.3}b) tells
us that $\delta_{0,b}:A \ra \langle b \rangle$ has infinite degree if
$\delta_{0,1}=\mu\circ\delta_{0,b}:A \ra Z_q$ has infinite degree. So, we only need
to show that $\delta_{0,1}:A \ra Z_q$
has infinite degree, if $A$ is infinite.\smallskip\\
\emph{Case 1, $A$ has infinite exponent:} In this case, $A$ contains an element
$a_n$ of order greater than $n$, for every $n\in\Z^+\!$. It follows that the elements
$0a_n,1a_n,\dotsc,na_n$ are pairwise different, so that
\[ [\Delta_{a_n}^n\delta_{0,b}](0) = \sum_{j=0}^n (-1)^{n-j} \binom{n}{j} \delta_{0,b}(0+ja_n)
=(-1)^n\delta_{0,b}(0)=(-1)^nb\neq0, \] as in Lemma \ref{0.0}.
Hence, $\fdeg(\delta_{0,1}) = \infty$, indeed. \smallskip\\
\emph{Case 2, $A$ is infinite and has finite exponent:} In this case, $\rank(A)=\infty$. So, for every
$n\in\Z^+$ there exists a subgroup $\underline{A}$ of $A$ with
$\rank(\underline{A})=n$. Applying first Proposition \ref{DELTAPROP}b), then
Theorem \ref{1.5} and then Lemma \ref{1.2}d), we see that
\[ \fdeg(\delta_{0,1}) = \delta(\underline{A},Z_p) = \nu(Z_p[\underline{A}])-1 \geq \rank(\underline{A})=n. \]
Hence, $\fdeg(\delta_{0,1}) = \infty$, as desired. Alternatively, this also follows from
the insight that, if $a_1, a_2, \dotsc, a_n$ generate $\underline{A}$, then the sum
$a_I:=\sum_{i\in I}a_i$ over a subset $I\subset \{1,2,\dotsc,n\}$ is zero only if
$I=\emptyset$, so that
$$[\Delta_{a_n}\Delta_{a_{n-1}}\dotsb\Delta_{a_1}\delta_{0,1}](0)f
=\sum_{I\subset \{1,\dotsc,n\}}\!\!(-1)^{n-|I|}\delta_{a_I,1}(0) =(-1)^n\delta_{0,1}(0)
=(-1)^n\neq0.
$$
b) We may assume that $A$ is finite, as otherwise part a) applies. So, there exists
an element $a\in A$ of order $p\in\mathcal{P}$ such that $B$ is not a $p$-group.
Corollary \ref{1.4}a) reduces us now further to the case $A = \langle a
\rangle\cong Z_p$. Moreover, $B$ contains an element $b$ of order $q$ coprime to
$p$. To show than $\fdeg(\delta_{0,b})=\infty$, we may also assume $B= Z_q$ and
$b=1\in Z_q$, exactly as in the proof of part a). In this setting, $\Hom(A,B) = \{0\}$
so that Corollary \ref{0.1.7} implies $\fdeg(\delta_{0,1}) = \infty$, as desired.\smallskip\\
c) First applying Theorem \ref{SUMTHM} with $A_i = Z_{p_i^{\alpha_i}}$, for
$i=1,\dotsc,r$, and then applying Theorem \ref{WEISMANWILSON}, we get
$$
\delta(A,B)
%\delta(\bigoplus_{i=1}^r Z_{p_i^{\alpha_i}},B) &
= \max_{\underline{\beta}} \sum_{i=1}^r\delta(Z_{p^{\alpha_i}},Z_{p^{\beta_i+1}})
= \sum_{i=1}^r(p^{\alpha_i}\!-1) + \max_{\underline{\beta}}
\sum_{i=1}^r \beta_i(p-1)p^{\alpha_i-1},
$$
where the maximum ranges over all $\underline{\beta} \in \N^r$ with
$|\underline{\beta}| \coloneqq \beta_1 + \ldots + \beta_r = \beta-1$.  However, for all
$\underline{\beta} \in \N^r$ with $|\underline{\beta}| = \beta-1$,
\[ \sum_{i=1}^r \beta_i(p-1)p^{\alpha_i-1} \leq \sum_{i=1}^r \beta_i(p-1)p^{\alpha_1-1}
= (\beta-1)(p-1)p^{\alpha_1-1}, \] with equality for $\underline{\beta} = (\beta-1,0,\ldots,0)$.
So, $\delta(A,B) = \sum_{i=1}^r (p^{\alpha_i}-1) + (\beta-1)(p-1)p^{\alpha_1-1}$.\smallskip\\
d) Let $f \in B^A$.  Since $A$ is finite and $B$ is a torsion group, the subgroup
$\underline{B}$ generated by $f(A)$ is finite.  By part c) the function
$f|^{\underline{B}}: A\ra\underline{B}$ has finite degree, so that by Corollary
\ref{1.4}b), we have $\fdeg(f) = \fdeg(f|^{\underline{B}})<\infty$. Finally, since $B$ is a
$p$-group of infinite exponent, $\exp(B[p^\beta]) = p^\beta$ for all $\beta \in \Z^+$.
Thus, Theorem\,\ref{MAINTHM1}\,c) implies
\[ \delta(A,B) \geq \smash{\sup_{\beta \in \Z^+}} \fdeg(A,B[p^\beta]) = \infty. \]
\end{proof}

\subsection{Computation of $\mathcal{D}(A,B)$} In this section we will compute $\mathcal{D}(A,B)$ for a class of commutative groups including all finitely
generated groups $A$ and $B$.  By Theorem \ref{MAINTHM1} and Lemma
\ref{0.1.5}, to compute $\mathcal{D}(A,B)$ for any $A$ and $B$ it remains to
determine $\delta^{\circ}(A,B)$.

\begin{prop}
\label{3.9} Let $A$ and $B$ be nontrivial commutative groups.
\begin{itemize}
\item[a)] If $A$ is a torsion group and $B$ is torsion free, then $\mathcal{D}(A,B) =
    \{-\infty,0,\infty\}$.\smallskip
\item[b)] If $A$ is a torsion group, $B$ is torsion-split with $B[\tors]\neq \{0\}$ and
    $B/B[\tors] \neq \{0\}$, then \[\mathcal{D}(A,B) = \mathcal{D}(A,B[\tors]) \cup
    \{\infty\}.\]
\item[c)] If there is a surjective homomorphism $\varepsilon: A \ra \Z$, then
    $\mathcal{D}(A,B) = \NN$.
\end{itemize}
\end{prop}
\begin{proof}
a) Since $\Hom(A,B) = \{0\}$, this follows from Corollary \ref{0.1.7}.\smallskip\\
b) We may write $B = B[\tors] \times B'$.  Let $\pi_1: B \ra B[\tors]$ and $\pi_2: B \ra
B'$ be the two projection maps. Then for $f\in B^A$, Lemma \ref{AM3.4} tells us that
$$\fdeg(f) = \max(\fdeg(\pi_1 \circ f), \fdeg(\pi_2 \circ f)).
$$
With this equation it is easy to verify that $\mathcal{D}(A,B) = \mathcal{D}(A,B[\tors])
\cup \{\infty\}$ by considering all possible degrees of functions $f_1: A \ra B[\tors]$
and $f_2: A \ra B'$, independently. By part a), $\mathcal{D}(A,B') = \{-\infty,0,\infty\}$,
as $B' \neq \{0\}$. So, $\mathcal{D}(A,B')\cap\N = \{0\}$ and this already shows that
$\mathcal{D}(A,B)\cap\N = \mathcal{D}(A,B[\tors])\cap\N$, as
$\mathcal{D}(A,B[\tors])\cap\N\neq\emptyset$. That
$\mathcal{D}(A,B)\cap\{-\infty,+\infty\} =
\mathcal{D}(A,B[\tors])\cap\{-\infty,+\infty\} \cup \{\infty\}$ follows from $\infty\in\mathcal{D}(A,B')$.\smallskip\\
c) Evidently $-\infty,0 \in \mathcal{D}(A,B)$, and by Theorem \ref{MAINTHM1}a) we
have $\infty \in \mathcal{D}(A,B)$.  Let $d \in \Z^+\!$.  Consider first the case
$A= B = \Z$.   We claim that the map $x\mapsto\tbinom{x}{d}$ from $\Z$ to $\Z$ %U several changes
has functional degree $d$. Indeed, by ``Pascal's Identity'' we have, for all $x
\in \Z$, that
\[ \Delta_1 \tbinom{x}{d} = \tbinom{x}{d-1}. \]
Hence, for a fixed $b_0 \in B \setminus \{0\}$, the map $g_d: \Z \ra B$ given by $x
\mapsto \tbinom{x}{d} b_0$ has functional degree $d$. So, by Lemma \ref{1.3}a), the
map $\varepsilon^* g: A \ra B$ also has functional degree $d$.
\end{proof}
\noindent Proposition \ref{3.9} reduces the computation of $\mathcal{D}(A,B)$ for
$A$ and $B$ finitely generated to the case in which $A$ and $B$ are finite.   The
next result accomplishes this for all torsion groups $A$ and $B$ such that
$\exp(A[p^{\infty}])$ is finite for all $p \in \mathcal{P}$, hence in particular for all finite
groups $A$.

\begin{thm}
\label{3.10} Let $A$ and $B$ be nontrivial torsion groups such that
$\exp(A[p^{\infty}])$ is finite for all $p \in \mathcal{P}$.
\begin{itemize}
\item[a)] We have
\[ \delta^{\circ}(A,B) = \sup_{p \in \mathcal{P}} \delta^{\circ}(A[p^{\infty}],B[p^{\infty}]). \]
\item[b)] If there are infinitely many $p \in \mathcal{P}$ such that $A[p^{\infty}]$
    and $B[p^{\infty}]$ are both nonzero, then
\[ \delta^{\circ}(A,B) = \infty. \]
\item[c)] If there is $p \in \mathcal{P}$ such that $A[p^{\infty}] \neq 0$ and
    $\exp(B[p^{\infty}]) = \infty$, then $\delta^{\circ}(A,B) = \infty$.\smallskip
\item[d)] If there is $p \in \mathcal{P}$ such that $A[p^{\infty}]$ is infinite and
    $B[p^{\infty}] \neq 0$, then $\delta^{\circ}(A,B) =\infty$.\smallskip
\item[e)] In the remaining case -- that is, the set of $p \in \mathcal{P}$ such that
    $A[p^{\infty}]$ and $B[p^{\infty}]$ are both nonzero is finite and for each such
    $p$, $A[p^{\infty}]$ is finite and $B[p^{\infty}]$ has finite exponent -- the
    parameter $\delta^{\circ}(A,B) = \sup_{p \in \mathcal{P}}
    \delta^{\circ}(A[p^{\infty}],B[p^{\infty}])$ is finite, and it can be computed using
    Theorem \ref{MAINTHM1}c).
\end{itemize}
\end{thm}
\begin{proof}
a) This follows from Theorem \ref{DIAGONALTHM}.\smallskip\\
b) If for $p \in \mathcal{P}$ we have that $A[p^{\infty}]$ and $B[p^{\infty}]$ are both
nonzero, then there is a surjective group homomorphism $\varepsilon: A \ra Z_p$
(here we use that $A[p^{\infty}]$ has finite exponent) and an injective group
homomorphism $\mu: Z_p \hookrightarrow B$, so by Theorem \ref{3.8}c), we
have \[\{-\infty,0,\ldots,p-1\} = \mathcal{D}(Z_p,Z_p) \subset \mathcal{D}(A,B). \]
Since this holds for infinitely
many $p \in \mathcal{P}$, we get $\delta^{\circ}(A,B) = \infty$.\smallskip\\
c) There is a surjective group homomorphism $\varepsilon: A \ra Z_p$ and for all
$\beta \in \Z^+$ an injective group homomorphism $\iota: Z_{p^\beta} \hookrightarrow B$,
so as above $\delta^{\circ}(A,B) \geq \delta^{\circ}(Z_p,Z_{p^\beta})$.  By Theorem
\ref{3.8}c), we have $\sup_\beta \delta^{\circ}(Z_p,Z_{p^\beta}) = \infty$, so
$\delta^{\circ}(A,B) = \infty$.\smallskip\\
d) It follows from Theorem \ref{PRUFERBAERTHM} that if $A[p^{\infty}]$ is infinite
and of bounded exponent then for each $d \in \Z^+$ there is a surjective
homomorphism $\varepsilon: A \ra \bigoplus_{i=1}^d Z_p$ and an injective group
homomorphism $\mu: Z_p \hookrightarrow B$.  Using Theorem \ref{3.8}c) as
above, it follows that
\[  \delta^{\circ}(A,B) \geq \sup_d \delta \left(\bigoplus_{i=1}^d Z_p,Z_p\right) = \infty. \]
e) This follows from Theorem \ref{3.8}c) and part a).
\end{proof}

\begin{example}
For $p \in \mathcal{P}$, let $C_{p^{\infty}} = \C^{\times}[p^{\infty}]$ be the
\textbf{Pr\"ufer $p$-group}, a $p$-group of infinite exponent. As the identity map is a nonzero
homomorphism, we have $\delta^{\circ}(C_{p^{\infty}},C_{p^{\infty}}) \geq 1$, but we do not
know more. This explains the need for the hypothesis that $\exp(A[p^{\infty}])$ is
finite for all $p \in \mathcal{P}$ in Theorem \ref{3.10}.
\end{example}

\end{document}